\documentclass[a4paper,12pt]{article}

\usepackage{amscd,amssymb,amsmath,amsfonts,amsthm, mathrsfs}
\usepackage{stmaryrd, accents}
\usepackage{graphicx}
\usepackage{verbatim,enumerate,paralist}
\usepackage{textcomp}
\usepackage{tensor}
\usepackage[colorlinks]{hyperref}
\hypersetup{%
  urlcolor=blue,linkcolor=blue,citecolor=blue
}
\usepackage{pdfpages}

\usepackage[T1]{fontenc}

\usepackage[scr=rsfs,cal=boondox]{mathalfa}

\usepackage[arrow, matrix, tips, curve, graph, rotate]{xy}
\SelectTips{cm}{10}

\makeatletter
\def\matrixobject@{%
  \edef \next@{={\DirectionfromtheDirection@ }}%
  \expandafter \toks@ \next@ \plainxy@
  \let\xy@@ix@=\xyq@@toksix@
  \xyFN@ \OBJECT@}
\let\xy@entry@@norm=\entry@@norm
\def\entry@@norm@patched{%
  \let\object@=\matrixobject@
  \xy@entry@@norm }
\AtBeginDocument{\let\entry@@norm\entry@@norm@patched}
\makeatother

\newcommand{\twocong}[2][0.5]{\ar@{}[#2] \save ?(#1)*{\cong}\restore}
\newcommand{\twoeq}[2][0.5]{\ar@{}[#2] \save ?(#1)*{=}\restore}
\newcommand{\ltwocell}[3][0.5]{\ar@{}[#2] \ar@{=>}?(#1)+/r 0.2cm/;?(#1)+/l 0.2cm/^{#3}}
\newcommand{\rtwocell}[3][0.5]{\ar@{}[#2] \ar@{=>}?(#1)+/l 0.2cm/;?(#1)+/r 0.2cm/^{#3}}
\newcommand{\utwocell}[3][0.5]{\ar@{}[#2] \ar@{=>}?(#1)+/d  0.2cm/;?(#1)+/u 0.2cm/_{#3}}
\newcommand{\dtwocell}[3][0.5]{\ar@{}[#2] \ar@{=>}?(#1)+/u  0.2cm/;?(#1)+/d 0.2cm/^{#3}}
\newcommand{\ultwocell}[3][0.5]{\ar@{}[#2] \ar@{=>}?(#1)+/dr  0.2cm/;?(#1)+/ul 0.2cm/^{#3}}
\newcommand{\urtwocell}[3][0.5]{\ar@{}[#2] \ar@{=>}?(#1)+/dl  0.2cm/;?(#1)+/ur 0.2cm/^{#3}}
\newcommand{\dltwocell}[3][0.5]{\ar@{}[#2] \ar@{=>}?(#1)+/ur  0.2cm/;?(#1)+/dl 0.2cm/^{#3}}
\newcommand{\drtwocell}[3][0.5]{\ar@{}[#2] \ar@{=>}?(#1)+/ul  0.2cm/;?(#1)+/dr 0.2cm/^{#3}}

\newcommand{\pullbackcorner}[1][dr]{\save*!/#1-1.2pc/#1:(-1,1)@^{|-}\restore}

\makeatletter
\newcommand{\xRightarrow}[2][]{\ext@arrow 0359\Rightarrowfill@{#1}{#2}}
\makeatother

\newtheorem{theorem}{Theorem}[section]
\newtheorem{corollary}[theorem]{Corollary}
\newtheorem{lemma}[theorem]{Lemma}
\newtheorem{proposition}[theorem]{Proposition}
\newtheorem{remark}[theorem]{Remark}
\newtheorem{definition}[theorem]{Definition}
\newtheorem{example}[theorem]{Example}

\newtheoremstyle{step}{2\bigskipamount}{\medskipamount}{\upshape}{}{\itshape}{. }{ }{\underline{Step~\thestep}}
\theoremstyle{step}

\renewcommand{\thestep}{\arabic{step}}

\newcommand{\lra}{\longrightarrow}

\newcommand{\Ra}{\Rightarrow}

\newcommand{\ldual}[1]{\mathord{{\let\nolimits\relax\sideset{^\wedge}{}{#1}}}}
\newcommand{\laction}[2]{\mathord{{\let\nolimits\relax\sideset{^{#1}}{}{#2}}}}
\newcommand{\conj}[2]{\mathord{{\let\nolimits\relax\sideset{^{#1}}{}{#2}}}}

\newcommand{\xra}{\xrightarrow}
\newcommand{\xla}{\xleftarrow}
\makeatletter
\newcommand{\xRa}[2][]{\ext@arrow 0359\Rightarrowfill@{#1}{#2}}
\makeatother

\usepackage[utf8]{inputenc}
\DeclareFontFamily{U}{min}{}
\DeclareFontShape{U}{min}{m}{n}{<-> udmj30}{}

\def\CA{{\mathscr A}}
\def\CB{{\mathscr B}}
\def\CC{{\mathscr C}}
\def\CD{{\mathscr D}}
\def\CE{{\mathscr E}}

\def\CK{{\mathscr K}}
\def\CL{{\mathscr L}}

\def\CO{{\mathscr O}}
\def\CP{{\mathscr P}}

\def\CR{{\mathscr R}}

\def\CT{{\mathscr T}}

\def\CX{{\mathscr X}}

\DeclareMathAlphabet{\mathbbe}{U}{bbold}{m}{n}

\begin{document}

\author{Ross Street}

\title{Objective Mackey and Tambara functors via parametrized categories}
\date{\today}
\maketitle

\tableofcontents
\newpage
\begin{abstract}
\noindent The first word in the title is intended in a sense suggested by
Lawvere and Schanuel whereby finite sets are objective natural numbers. 
At the objective level, the axioms defining abstract Mackey and Tambara 
functors are categorically familiar. The first step was taken by Harald
Lindner in 1976 when he recognized that Mackey functors, defined as
pairs of functors, were equivalently single functors with domain a category of spans.
In 1993 Tambara recognized that TNR-functors (that is, functors designed to
have abstract trace, norm and restriction operations, and now called Tambara functors)
were equivalently certain functors out of a category of polynomials.
We define objective Mackey and objective Tambara functors as parametrized 
categories which have local finite products and satisfy some parametrized 
completeness and cocompleteness restriction. However, we can replace 
the original parametrizing base for objective Mackey functors by a 
bicategory of spans while the replacement for objective Tambara functors 
is a bicategory obtained by iterating the span construction; 
these iterated spans are polynomials. There is an objective Mackey functor 
of ordinary Mackey functors. We show that there is a distributive law 
relating objective Mackey functors to objective Tambara functors analogous 
to the distributive law relating abelian groups to commutative rings. 
We remark on hom enrichment matters involving the 2-category 
$\mathrm{Cat}_{+}$ of categories admitting finite coproducts and functors 
preserving them, both as a closed base and as a skew-closed base.
\end{abstract}

\section*{Introduction}

The sense in which we use the terms ``abstract'' and ``objective'' is that of Lawvere and Schanuel (for example, see \cite{Schanuel}). Coproduct and product of finite sets have
the virtue, over abstract addition and multiplication of natural numbers, of universal
properties which can help in understanding of combinatorial equations as bijections.   

The process of moving from objective to abstract is very old. 
Like quantization in Physics, the reverse process is
also rather old, yet it is not well defined even up to any reasonable notion of equivalence. 
The term ``categorification'' was coined by Crane--Frenkel in low-dimensional topology (see \cite{CraFre}) and then adopted by representation theorists.  

The Introduction of David Chan's paper \cite{Ch} provides good background to the roles of
Mackey and Tambara functors in equivariant stable homotopy theory.
In particular, the third paragraph begins with the following sentence:
\newline {\em While Mackey functors should be thought of as the correct analogue of an
abelian group in the equivariant setting, the role of commutative rings is played
by objects known as Tambara functors.} 


Motivated by terminology and constructions in \cite{HH, BH18, BH22, Ch}, my goal was to 
investigate in what way Tambara functors on a lextensive category $\CE$ are 
Mackey functors with extra structure. However, I soon moved to the objective level (in a different direction from \cite{BD'A})
and found what I consider is a satisfactory confirmation 
of the analogy between Mackey functors and abelian groups, and commutative rings and Tambara functors.
For the case where $\CE$ is a category of sets on which a group acts, I see that a left adjoint to 
an inclusion of Mackey functors into Tambara functors, thought of as the analogue of the monoid-ring functor, 
was constructed and applied by Nakaoka in \cite{Nak}.  
 
The main idea is that the basic construction is a variant of 
B\'enabou's bicategory $\mathrm{Spn}\CE$ of spans $U\xla{u}S\xra{v}V$ in a category $\CE$ 
with pullbacks (see \cite{Ben1967}).
We allow $\CE$ to be replaced by a bicategory $\CC$ equipped with a chosen wide 
(= all object inclusive), locally full subbicategory $\CR$ with a few conditions making it what we here call a {\em protocalibration}.
For our bicategory $\mathrm{Spn}_{\CR}\CC$, we restrict our 1-morphisms to spans for which the 
left leg $u$ is in $\CR$ yet we allow our 2-morphisms to be isomorphism classes of span morphisms which are pseudo on the left and lax on the right. Its universal property is that categories parametrized by
$\mathrm{Spn}_{\CR}\CC$ are essentially categories parametrized by $\CC$ which are 
$\CR$-cocomplete. This last term refers to the existence of pointwise left extensions along 
$\CC$-parametrized functors of the form $\CC(-,U)\xra{\CC(-,r)}\CC(-,V)$ with $r\in \CR$.   

Other constructions are obtained by duality and iteration. If $\CL$ is a proto-calibration of
$\CC^{\mathrm{co}}$ then categories parametrized by $\tensor[_{\CL}]{\mathrm{Spn}}{}\CC = (\mathrm{Spn}_{\CL}\CC^{\mathrm{co}})^{\mathrm{co}}$ are essentially categories 
parametrized by $\CC$ which are $\CL$-complete. Moreover, with the specific iteration $\mathrm{Spn}_{\CE_*}(\tensor[_{\CP}]{\mathrm{Spn}}{}\CE)$ we obtain the opposite of the bicategory
$\mathrm{Ply}\CE$ of polynomials in a category $\CE$ with pullbacks where $\CP$ consists of the
powerful morphisms in $\CE$ and $\CE_*$ consists of spans with the left leg $u$ invertible.
The universal property (in the case where $\CE$ is locally cartesian closed) obtained by iteration
agrees with that proved by Charles Walker \cite{Walker2019}.

Rings are Eilenberg-Moore algebras for a monad on $\mathrm{Set}$ obtained as a composite monad
using a Beck distributive law \cite{DL} of the monad for monoids over the monad for abelian groups. 
In Proposition~\ref{monadicityofobjTambara} we see that objective Tambara functors are Eilenberg-Moore pseudo-algebras 
for a pseudomonad on the bicategory of finite (bicategorical) product preserving pseudofunctors 
$\CE^{\mathrm{op}}\to \mathrm{Cat}$ obtained as a composite pseudomonad using a 
pseudodistributive law over the pseudomonad for objective Mackey functors. 

\section{Parametrized categories and protocalibrations}\label{Ropc}

Categories varying over \cite{17} (fibred over \cite{Ben1975b}, indexed by \cite{ICaTA}, or parametrized by \cite{32}) a category or bicategory $\CC$
can be regarded as pseudofunctors $\mathbb{X} : \CC^{\mathrm{op}}\to \mathrm{Cat}$.
We think of an object of $\mathbb{X}U$ as a {\em generalised object} of the $\CC$-variable category $\mathbb{X}$.
We write 
$$\mathrm{CAT}\CC = \mathrm{Ps}( \CC^{\mathrm{op}}, \mathrm{Cat})$$ 
for the (strict) bicategory of pseudofunctors
(called ``homomorphisms'' in \cite{Ben1967}), pseudonatural transformations, and modifications. 

A $\CC$-variable category $\mathbb{X}\in \mathrm{CAT}\CC$ is {\em cocomplete with respect to a morphism}
$j : \mathbb{A} \to \mathbb{B}$ in $\mathrm{CAT}\CC$ (see Section 7 of \cite{17}) when every morphism $f : \mathbb{A} \to \mathbb{X}$
has a pointwise (see Section 4 of \cite{8}) left extension along $j$. 
Generally, pointwiseness translates to a lax form of the Chevalley-Beck (CB) condition. However, the pseudo form is relevant when dealing with groupoid fibrations as we will now remind the reader. 

A morphism $\kappa : z\to x$ in a category $E$ is {\em cartesian} for a functor $p : E \to B$ when the following square \eqref{cartmor} is a pullback in $\mathrm{Set}$ for all objects $k$ of $E$.
\begin{equation}\label{cartmor}
 \begin{aligned}
\xymatrix{
E(k,z) \ar[d]_{p}^(0.5){\phantom{aaaaaaaa}}="1" \ar[rr]^{E(k,\kappa)}  && E(k,x) \ar[d]^{p}_(0.5){\phantom{aaaaaaaa}}="2"
\\
B(pk,pz) \ar[rr]_-{B(pk,p\kappa)} && B(pk,px) 
}
 \end{aligned}
\end{equation} 
A functor $p : E\to B$ is a {\em groupoid fibration (over $B$)} when every morphism of $E$ is cartesian for $p$ and, for all $x\in E$ and $\phi : b\to px$ in $B$, there
exists a morphism $\kappa : z\to x$ in $E$ and an isomorphism $b\cong pz$ whose
composite with $p\kappa$ is $\phi$.
A functor $p : E\to B$ is a {\em groupoid opfibration (over $B$)} when $p : E^{\mathrm{op}}\to B^{\mathrm{op}}$ is a groupoid fibration.

A morphism $p : E\to B$ in a bicategory $\CC$ is a {\em groupoid fibration} when, for all objects $U\in \CC$, the functor $\CC(U,p) : \CC(U,E)\to \CC(U,B)$ is a groupoid fibration. The proof of Proposition~\ref{ptwsalonggpdfib} is essentially as for Proposition 23 of \cite{8}; alternatively,
check it in $\mathrm{Cat}$ and use Yoneda since all concepts are representable.  

\begin{proposition}\label{ptwsalonggpdfib} 
In any finitely complete bicategory, suppose $p : E\to B$ is a groupoid opfibration.
In diagram~\eqref{pseudoptwise}, the 2-morphism $\sigma$ exhibits $k$ as a pointwise left extension
of $f$ along $p$ if and only if, for all morphisms $b$ with the square a bicategorical pullback (= briefly, bipullback), the pasted composite at $p$ exhibits $kb$ as a left extension of $f\tilde{b}$
along $\tilde{p}$.  
\begin{equation}\label{pseudoptwise}
 \begin{aligned}
\xymatrix{
E \ar[rd]_{f}^(0.5){\phantom{a}}="1" \ar[rr]^{p}  && B \ar[ld]^{k}_(0.5){\phantom{a}}="2" \ar@{=>}"1";"2"^-{\sigma}
\\
& X 
}
\qquad
\xymatrix{
P \ar[d]_{\tilde{b}}^(0.5){\phantom{aaaa}}="1" \ar[rr]^{\tilde{p}} \pullbackcorner  && V \ar[d]^{b}_(0.5){\phantom{aaaa}}="2" \ar@{=>}"1";"2"^-{ }_-{\cong}
\\
E \ar[rr]_-{p} && B 
}
 \end{aligned}
\end{equation}    
\end{proposition}

If $\CR$ is a set of morphisms of a bicategory $\CC$, we say $\mathbb{X}$ is {\em $\CR$-cocomplete} when it is cocomplete with respect to all morphisms 
$\CC(-,j) : \CC(-,U)\to \CC(-,V)$ for $j : U\to V$ in $\CR$.
A pseudonatural transformation  $\theta : \mathbb{X}\to \mathbb{Y}$ is {\em $\CR$-cocontinuous} when it preserves left extensions along the $\CC(-,j)$ for $j\in \CR$. 

Let $\mathrm{CAT}_{\CR}\CC$ denote the subbicategory of $\mathrm{CAT}\CC$ consisting of the $\CR$-cocomplete objects and $\CR$-cocontinuous morphisms.

We are interested in those $\CR$ in possession of most of the properties in common with the ``calibrations'' of \cite{134}, and the calibrations of B\'enabou \cite{Ben1975a} when $\CC$ is a category.
In the case where $\CC$ is a category of $G$-sets, these conditions plus closure under coproducts (which we will also require in Section~\ref{Fcfsap}) are those called ``indexing categories'' 
in equivariant stable homotopy theory; see Chan \cite{Ch}.  

\begin{definition}\label{protocalib}
A {\em protocalibration} of a bicategory $\CC$ is a set $\CR$ of morphisms when 
it satisfies the four conditions:
\begin{itemize}
\item[I.] all equivalences are in $\CR$ and any morphism isomorphic to one in $\CR$ is in $\CR$;
\item[C.] composites of morphisms in $\CR$ are in $\CR$;
\item[P.] for all $U\xra{f} W$ in $\CR$ and all $V\xra{g} W$ in $\CC$, the bicategorical pullback 
$P\xra{\bar{f}} V$ of $f$ along $g$ exists and is in $\CR$;
\item[G.] every morphism in $\CR$ is a groupoid opfibration.
\end{itemize} 
Each protocalibration $\CR$ will be identified with the subbicategory
of $\CC$ containing all objects, only the morphisms in $\CR$, and all 2-morphisms between them.
Notice that a protocalibration of $\CC^{\mathrm{co}}$ is the same as one for $\CC$ except for 
condition (G) where ``groupoid opfibration'' must be replaced by ``groupoid fibration''.  
\end{definition}

\begin{example}\label{protocalibexx}
\begin{itemize}
\item[1.] Suppose $\CC$ is just a category (locally discrete as a bicategory).
Then every morphism is both a groupoid fibration and a groupoid opfibration, 
and bicategorical pullbacks are pullbacks.
If $\CC$ has pullbacks, the set of all morphisms is a protocalibration.
\item[2.] Again suppose $\CC$ is a category. If pullbacks of powerful (= exponentiable)
morphisms exist then the powerful morphisms form a protocalibration (see Corollary 2.6 of \cite{104}).
This applies to the {\em category} $\CC=\mathrm{Cat}$ of categories and functors where the powerful morphisms are the Giraud-Conduch\'e functors \cite{FCond}.
\item[3.] Let $\CE$ be a category with pullbacks and every morphism powerful. The groupoid
fibrations (= left-adjoint morphisms) in the bicategory $\mathrm{Spn}\CE$ of spans in $\CE$ is
a protocalibration of $(\mathrm{Spn}\CE)^{\mathrm{co}}$ (see Propositions 4.1 and 8.3 of \cite{134}). 
\end{itemize}
\end{example}

\begin{proposition}[Compare \cite{13} Proposition 9.12, \cite{17} Proposition 7.12]\label{CB}  
Let $\CR$ be a protocalibration of the bicategory $\CC$. A $\CC$-variable category $\mathbb{X}$ is $\CR$-cocomplete if and only if, for all $f \in \CR$, the functor $\mathbb{X}f$ has a left adjoint $\mathbb{X}_*f$ and, for all
bicategorical pullback squares
\begin{equation}\label{bipb}
 \begin{aligned}
\xymatrix{
P \ar[d]_{g_f}^(0.5){\phantom{aaaa}}="1" \ar[rr]^{f_g} \pullbackcorner  && V \ar[d]^{g}_(0.5){\phantom{aaaa}}="2" \ar@{=>}"1";"2"^-{ }_-{\cong}
\\
U \ar[rr]_-{f} && W 
}
 \end{aligned}
\end{equation}  
in $\CE$, the mate $$\mathbb{X}_*f_g \ \circ \ \mathbb{X}g_f\ \Ra \ \mathbb{X}g \ \circ \ \mathbb{X}_*f $$ of the isomorphism
$$\mathbb{X}g_f \ \circ \  \mathbb{X}f \ \cong \ \mathbb{X}f_g \ \circ \ \mathbb{X}g$$ is invertible. A pseudonatural transformation $\theta : \mathbb{X}\to \mathbb{Y}$ is $\CR$-cocontinuous if and only if, for all $U\xra{f}W$ in $\CR$, the mate $\mathring{\theta}_f : \mathbb{Y}_*f\circ \theta_U\to \theta_W\circ \mathbb{X}_*f$of $\theta_f : \theta_U\circ \mathbb{X}f\to \mathbb{Y}f\circ \theta_W$ is invertible.   
\end{proposition}
\begin{remark}\label{lazy} 
Henceforth we occasionally allow ourselves to omit the invertible 2-morphisms in bipullbacks such as \eqref{bipb}.
\end{remark}

\begin{lemma}\label{finprodofmorphsinL}
Suppose $\CC$ has binary bicategorical products\footnote{We will abbreviate this to ``biproduct'' and use ``direct sum'' when product and coproduct coincide.} and $\CR$ is a protocalibration of $\CC$.
If $X\xra{f}Y$ and $X'\xra{f'}Y'$ are in $\CR$ then so too is $X\times X'\xra{f\times f'}Y\times Y'$. So, if $\CC$ has finite bicategorical products then $\CR$ is a monoidal subbicategory of the cartesian monoidal bicategory $\CC$. 
\end{lemma}
\begin{proof} The morphism $X\times X'\xra{f\times f'}Y\times Y'$ is the composite
of $f\times 1_{X'}$ and $1_Y\times f'$ which are respectively bipullbacks of $f$ and $f'$
along projections. 
\end{proof}

\begin{corollary}\label{Dress}
Suppose $\CC$ has binary biproducts and $\CR$ is a protocalibration of $\CC$.
If $\mathbb{Z}$ is $\CR$-cocomplete then so is $\mathbb{Z}(-\times W)$ for all $W\in \CC$. 
\end{corollary}
\begin{proof} Cocompleteness of $\mathbb{Z}(-\times W)$ with respect to any $f$ amounts, by Yoneda, to cocompleteness of $\mathbb{Z}$ with respect to $f\times 1_W$.  
\end{proof}

\begin{proposition}\label{inthom}
Suppose $\CC$ has binary biproducts and $\CR$ is a protocalibration of $\CC$.
For $\mathbb{Y}$ and $\mathbb{Z}$ in $\mathrm{CAT}_{\CR}\CC$, the object
$[\mathbb{Y},\mathbb{Z}]_{\CR} \in \mathrm{CAT}\CC$, defined on objects by
\begin{eqnarray*}
[\mathbb{Y},\mathbb{Z}]_{\CR}W = \mathrm{CAT}_{\CR}\CC(\mathbb{Y},\mathbb{Z}(-\times W)) 
\end{eqnarray*}
and on morphisms $W'\xra{h}W$ by composition with $\mathbb{Z}(-\times h)$, is in $\mathrm{CAT}_{\CR}\CC$.
\end{proposition}
\begin{proof} By the last sentence of Proposition~\ref{CB}, we need to see that, if $\theta : \mathbb{Y}\to \mathbb{Z}(-\times W)$
is $\CR$-cocontinuous, so too is $\mathbb{Z}(-\times h)\circ \theta$. Take $V\xra{r}U$
in $\CR$. The right square in the diagram  
\begin{eqnarray*}
 \begin{aligned}
\xymatrix{
\mathbb{Y}U \ar[d]_{\mathbb{Y}r}^(0.5){\phantom{aaaaaa}}="1" 
\ar[rr]^{\theta_U}  
&& \mathbb{Z}(U\times W) \ar[rr]^{\mathbb{Z}(1_U\times h)} 
     \ar[d]^{\mathbb{Z}(r\times 1_W)}_(0.5){\phantom{aaaaaa}}="2" \ar@{=>}"1";"2"^-{\theta_r}_-{\cong} 
     \ar@{}[d]^{\phantom{aaaaaaaaaa}}="3"
&& \mathbb{Z}(U\times W') \ar[d]^{\mathbb{Z}(r\times 1_{W'})}_(0.5){\phantom{aaaaaaaaaa}}="4" \ar@{=>}"3";"4"^-{}_-{\cong}
\\
\mathbb{Y}V \ar[rr]_-{\theta_V} && \mathbb{Z}(V\times W) \ar[rr]_-{\mathbb{Z}(1_V\times h)} && \mathbb{Z}(V\times W') 
}
 \end{aligned}
\end{eqnarray*}
is the value of $\mathbb{Z}$ on a bipullback in $\CC$. The invertibility of the mate of the
pasted 2-morphism follows from $\CR$-cocontinuity of $\theta$ and $\CR$-cocompleteness
of $\mathbb{Z}$. 
\end{proof}

\section{$\CR$-cocompletion}
Suppose $\CR$ is a protocalibration of the bicategory $\CC$. The {\em $\CR$-cocompletion $\overrightarrow{\mathbb{X}}$ of $\mathbb{X}\in \mathrm{CAT}\CC$} is described as follows (see \cite{4, Ben1975b, 17}). An object of the category $\overrightarrow{\mathbb{X}}U$ is a pair $(U\xla{u}S, x\in \mathbb{X}S)$ with $u\in \CR$; we might think of such an object as a $U$-indexed
family of generalised objects of $\mathbb{X}$.   
A morphism $$(U\xla{u}S, x\in \mathbb{X}S)\to (U\xla{u'}S', x'\in \mathbb{X}S')$$
 of $\overrightarrow{\mathbb{X}}U$ is an isomorphism class of triples $(S\xra{w}S', u'w\cong u, x\xra{\xi} (\mathbb{X}w)x')$
 where $\xi$ is a morphism of $\mathbb{X}S$; if we identify representables with their representing objects, the picture of such a morphism is a lax morphism \eqref{Rcocompletion} of spans in $\mathrm{CAT}\CC$. 
    \begin{eqnarray}\label{Rcocompletion}
\begin{aligned}
\xymatrix{
& & S \ar[d]_{ w} \ar[lld]_-{u\in \CR} \ar[rrd]^-{ x}  
\ar @{} [rd] | {\stackrel{\xi} \Leftarrow} & &
\\
U    && S' \ar[ll]^-{u'\in \CR} \ar[rr]_-{x'} \ar @{} [lu] | {\stackrel{\cong} \Leftarrow} & & \mathbb{X} }  
\end{aligned}
\end{eqnarray}
For $r : V\to U$ in $\CC$, define the functor $\overrightarrow{\mathbb{X}}r :\overrightarrow{\mathbb{X}}U\to \overrightarrow{\mathbb{X}}V$ to take $(u,x)$ to
$(u_r, (\mathbb{X}r_u)x)$ where
\begin{eqnarray}\label{cocopbs}
\begin{aligned}
\xymatrix{
P \ar[rr]^-{u_r} \ar[d]_-{r_u} \pullbackcorner && V \ar[d]^-{r} \\
S \ar[rr]_-{u} && U }
\qquad
\xymatrix{
P' \ar[rr]^-{u'_{r}} \ar[d]_-{r_{u'}} \pullbackcorner && V \ar[d]^-{r} \\
S' \ar[rr]_-{u'} && U}
\end{aligned}
\end{eqnarray}
are bipullbacks, and $(w,\xi)$ to $(\bar{w},\bar{\xi})$ where $\bar{w}$ is defined by compatible $u'_r\bar{w} \cong u_r$,
$r_{u'}\bar{w} \cong wr_u$, and $\bar{\xi}$ by 
$$\bar{\xi} : = \left( (\mathbb{X}r_u)x\xra{(\mathbb{X}r_u)\xi}(\mathbb{X}r_u)(\mathbb{X}w)x' \cong (\mathbb{X}\bar{w})(\mathbb{X}r_{u'})x'\right) \ .$$
For a 2-morphism $\alpha : r\Ra s : V\to U$, let $S\xla{s_u}Q\xra{u_s}V$ be the bipullback of the cospan $S\xra{u}U\xla{s}V$. Since $u$ is a groupoid opfibration, there exists a cocartesian 2-morphism $\kappa : r_u\Ra k$ and isomorphism $su_r\cong uk$ such that $u\kappa$ is the composite
$ur_u\cong ru_r\xra{\alpha u_r}su_r\cong uk$.  
If $r\in \CR$ then $\overrightarrow{\mathbb{X}}r$ has a left adjoint $\overrightarrow{\mathbb{X}}_*r : \overrightarrow{\mathbb{X}}V\to \overrightarrow{\mathbb{X}}U$ taking $(V\xla{v}T, y\in \mathbb{X}T)$ to
$(U\xla{rv}T, y\in \mathbb{X}T)$ (noting that $rv\in \CR$ by condition (C)).
To verify the pseudo-Chevalley-Beck condition, we refer to the bipullbacks \eqref{CBpbs} with 
$f, u\in \CR$. 
    \begin{eqnarray}\label{CBpbs}
\begin{aligned}
\xymatrix{
Q \ar[rr]^-{u_{g_f}} \ar[d]_-{g_{fu}} \pullbackcorner && P \ar[rr]^-{f_g} \ar[d]_-{g_f} \pullbackcorner && V \ar[d]^-{g} \\
S \ar[rr]_-{u\in \CR} && U  \ar[rr]_-{f\in \CR} && W}
\end{aligned}
\end{eqnarray}
Then we see that
\begin{eqnarray*}
(\overrightarrow{\mathbb{X}}g\circ \overrightarrow{\mathbb{X}}_*f) (u,x) & = & (\overrightarrow{\mathbb{X}}g) (fu,x) \  \cong \  (f_gu_{g_f}, (\mathbb{X}g_{fu})x) \\
& = & (\overrightarrow{\mathbb{X}}_*f_g)(u_{g_f},(\mathbb{X}g_{fu})x) \\ 
& \cong & (\overrightarrow{\mathbb{X}}_*f_g\circ \mathbb{X}g_f)(u,x)     
\end{eqnarray*}
naturally in $(u,x)\in \overrightarrow{\mathbb{X}}U$. So $\overrightarrow{\mathbb{X}}\in \mathrm{CAT}\CE$ is $\CR$-cocomplete by Proposition~\ref{CB}. 

\begin{remark}\label{ElikeSet} Suppose $\CC$ is a category with pullbacks so that we can take $\CR$ to contain
all morphisms. The terminal object $\mathbf{1}$ of $\mathrm{CAT}\CC$ takes all
objects of $\CC$ to the terminal category. 
Then $\mathbb{C} : = \overrightarrow{\mathbf{1}}$
takes $U\in \CC$ to the slice category $\CC_{/ U}$ and is defined on morphisms using pullback. Notice that an ordinary category has coproducts if and only if
all pointwise left Kan extensions into it exist along functors between discrete categories. Also $\mathrm{Set}$ is the coproduct cocompletion of the terminal ordinary category. Now we see that $\mathbb{C}$ is the coproduct cocompletion of the terminal $\CC$-variable category. This is one justification for viewing $\mathbb{C}$ as playing the role in $\CC$-variable category theory of $\mathrm{Set}$ in ordinary category theory. Also $\mathbb{C}$ is the internal full subcategory of $\CC$ at the heart of the
Yoneda structure for $\CC$-variable category theory; see \cite{13}. 
\end{remark}
\begin{definition}\label{C_R} For a protocalibration $\CR$ of a bicategory $\CC$, we write $\mathbb{C}_{\CR}$ for the $\CR$-cocompletion $\overrightarrow{\mathbf{1}}$ of $\mathbf{1}$. When $\CC$ has bipullbacks we write $\mathbb{C}$ for
$\mathbb{C}_{\CC}$, but notice that $\mathbb{C}_*$ exists even without the bipullbacks since it is defined on morphism using
composition.      
\end{definition}
\begin{remark}\label{laxcommaview}
The coproduct completion $\mathrm{Fam}X$ of a category $X$ can be constructed as the lax comma category
\begin{equation*}
 \begin{aligned}
  \xymatrix{
\mathrm{Fam}X \ar[d]_{}^(0.5){\phantom{aaaaa}}="1" \ar[rr]^{!}  && \mathbf{1} \ar[d]^{X}_(0.5){\phantom{aaaaa}}="2" \ar@{~>}"1";"2"^-{}
\\
\mathrm{Set} \ar[rr]_-{\mathrm{incl.}} && \mathrm{Cat} \ .
}
 \end{aligned}
\end{equation*}
Analogously, the $\CR$-cocompletion $\overrightarrow{\mathbb{X}}$ of $\mathbb{X}\in \mathrm{CAT}\CC$ can be constructed as the lax comma category
 \begin{equation*}
 \begin{aligned}
  \xymatrix{
\overrightarrow{\mathbb{X}} \ar[d]_{}^(0.5){\phantom{aaaaaa}}="1" \ar[rr]^{}  && \CC(-,\mathbf{1}) \ar[d]^{\mathbb{X}}_(0.5){\phantom{aaaaaa}}="2" \ar@{~>}"1";"2"^-{ }
\\
\mathbb{C}_{\CR} \ar[rr]_-{\mathrm{Yon}\restriction} && \mathrm{Ps}(\mathbb{C}_*^{\mathrm{op}},\mathrm{Cat})
}
 \end{aligned}
\end{equation*}
in $\mathrm{Ps}(\CC^{\mathrm{op}},\mathrm{CAT})$, where the bottom arrow is the restricted Yoneda embedding and the right arrow corresponds to $\mathbb{X}$ under the bicategorical Yoneda Lemma.  
\end{remark}

Directly from the explicit descriptions, we can verify the following.

\begin{proposition}\label{Gro_constr} For a protocalibration $\CR$ of a bicategory $\CC$
and a pseudofunctor $\mathbb{X} \in \mathrm{CAT}\CC$, the unique $! : \mathbb{X}\to \mathbf{1}$ induces a pseudonatural transformation $\mathbb{X}! : \overrightarrow{\mathbb{X}}\to \mathbb{C}_{\CR}$ which is a componentwise fibration. If $\CC$ is a finitely complete
category and each $! : S\to 1$ is in $\CR$ then $\mathbb{X}!_{1} : \overrightarrow{\mathbb{X}}1\to \mathbb{C}_{\CR}1$ is the fibration over $\CC$ constructed by Grothendieck from the pseudofunctor
$\mathbb{X} : \CC^{\mathrm{op}}\to \mathrm{Cat}$.     
\end{proposition}

There is a pseudonatural transformation $\eta_{\mathbb{X}} : \mathbb{X}\to \overrightarrow{\mathbb{X}}$ whose component $\eta_{\mathbb{X} U} : \mathbb{X}U\to \overrightarrow{\mathbb{X}}U$ at $U\in \CE$ is defined by $\eta_{\mathbb{X} U}x = (U\xla{1_U}U, x\in \mathbb{X} U)$.   

\begin{proposition}\label{ladjtounit} Suppose $\CR$ is a protocalibration of the bicategory $\CC$. Then
$\mathbb{X}\in \mathrm{CAT}\CC$ is $\CR$-cocomplete if and only if $\eta_{\mathbb{X}}$
has a left adjoint in $\mathrm{CAT}\CC$.   
\end{proposition}
\begin{proof} A left adjoint to $\eta_{\mathbb{X} U} : \mathbb{X}U\to \overrightarrow{\mathbb{X}}U$ precisely gives, for each $(u,x)\in \overrightarrow{\mathbb{X}}U$,
an object $u_*(x)\in \mathbb{X}U$ and an isomorphism 
$\mathbb{X}U(u_*(x), y) \cong \overrightarrow{\mathbb{X}}U((u,x),(1_U,y))$ natural in $y\in \mathbb{X}U$. However, 
$\overrightarrow{\mathbb{X}}U((u,x),(1_U,y))\cong \mathbb{X}S(x,(\mathbb{X}u)y)$.
So this amounts precisely to having a left adjoint $\mathbb{X}_*u$ to $\mathbb{X}u$
for all $u\in \CR$. Pseudonaturality of $\eta_{\mathbb{X}}$ amounts to the $\CR$-Chevalley-Beck condition for $\mathbb{X}$. The result follows now from Proposition~\ref{CB}.         
\end{proof}

We think of the left adjoint to $\eta_{\mathbb{X}}$ as assigning a coproduct to the indexed families of
generalised objects.

\begin{corollary}\label{KZness}
The pseudonatural transformation $\eta_{\overrightarrow{\mathbb{X}}} : \overrightarrow{\mathbb{X}}\to \overrightarrow{\overrightarrow{\mathbb{X}}}$ has a left adjoint
$\mu_{\mathbb{X}} : \overrightarrow{\overrightarrow{\mathbb{X}}}\to \overrightarrow{\mathbb{X}}$
in $\mathrm{CAT}\CC$.
\end{corollary}

Much as in Section 2 of \cite{8} and Section 3 of \cite{14}, for each protocalibration $\CR$ of $\CC$, 
we see we have a Kock-Z\"oberlein \cite{AKocka, AKockb, Zoeb} pseudomonad 
$$\overrightarrow{(-)}_{\CR} = (\overrightarrow{(-)}, \eta, \mu)$$ on $\mathrm{CAT}\CC$.

\begin{corollary}\label{pseudomonadicity}
The inclusion of the subbicategory $\mathrm{CAT}_{\CR}\CC$ of $\mathrm{CAT}\CC$, consisting of the $\CR$-cocomplete objects and $\CR$-cocontinuous morphisms, is pseudomonadic. That is, the inclusion is biequivalent over
$\mathrm{CAT}\CC$ to the forgetful pseudofunctor from the bicategory of
Eilenberg-Moore pseudo-algebras for the pseudomonad $\overrightarrow{(-)}_{\CR}$. 
In particular, $\mathrm{CAT}_{\CR}\CC$ is a complete bicategory and bicategorical limits
are preserved by the inclusion.   
\end{corollary}

The {\em opposite} $\mathbb{X}^{\mathrm{op}}$ of $\mathbb{X}\in \mathrm{CAT}\CC$ is obtained by composition with the 2-isomorphism $(-)^{\mathrm{op}} : \mathrm{Cat}\to \mathrm{Cat}^{\mathrm{co}}$. This is the value on objects for a 2-isomorphism
\begin{eqnarray}\label{opiso}
(-)^{\mathrm{op}} : \mathrm{CAT}\CC\to (\mathrm{CAT}\CC^{\mathrm{co}})^{\mathrm{co}} \ .
\end{eqnarray}
 
 Suppose $\CL$ is a protocalibration of $\CC^{\mathrm{co}}$.
 We define $\mathbb{X}\in \mathrm{CAT}\CC$ to be {\em $\CL$-complete} when
 $\mathbb{X}^{\mathrm{op}}$ is $\CL$-cocomplete. It follows that the 
 {\em $\CL$-completion}  $\overleftarrow{\mathbb{X}}$ of $\mathbb{X}\in \mathrm{CAT}\CC$ is 
$\overrightarrow{(\mathbb{X}^{\mathrm{op}})}^{\mathrm{op}}$.
If $\mathbb{X}$ is $\CL$-complete and $V\xra{f}U$ is in $\CL$, we write $\mathbb{X}_!f$ for the right adjoint of $\mathbb{X}f : \mathbb{X}U\to \mathbb{X}V$. 

We write $\tensor[_{\CL}]{\mathrm{CAT}\CC}{ }$ for the locally full subbicategory of $\mathrm{CAT}\CC$ consisting of the $\CL$-complete objects and $\CL$-continuous morphisms. The 2-isomorphism \eqref{opiso} induces an isomorphism
\begin{eqnarray}
\tensor[_{\CL}]{\mathrm{CAT}}{ }\CC \cong (\tensor[{}]{\mathrm{CAT}}{_{\CL} }\CC^{\mathrm{co}})^{\mathrm{co}} \ .
\end{eqnarray}

We have a dual Kock-Z\"oberlein pseudomonad on $\mathrm{CAT}\CC$ for which the Eilenberg-Moore pseudo-algebras are the $\CL$-complete objects where the algebra action
on the object is right adjoint to the unit. 

An object of $\overleftarrow{\mathbb{X}}U$ is a pair $(U\xla{u}S, x\in \mathbb{X}S)$ with $u\in \CL$, as for $\overrightarrow{\mathbb{X}}U$.   
A morphism $$(U\xla{u}S, x\in \mathbb{X}S)\to (U\xla{u'}S', x'\in \mathbb{X}S')$$
 of $\overleftarrow{\mathbb{X}}U$ is an isomorphism class of triples $(S'\xra{w}S, u'\cong uw, (\mathbb{X}w)x\xra{\zeta}x')$
 where $\zeta$ is a morphism of $\mathbb{X}S'$. 
    \begin{eqnarray}\label{Lcompletion}
\begin{aligned}
\xymatrix{
U \ar @{} [rrrd] | {\stackrel{\cong} \Rightarrow} & & S  \ar[ll]_-{u\in \CL} \ar[rr]^-{ x} \ar @{} [rd] | {\stackrel{\zeta} \Rightarrow }   & & \mathbb{X}
\\
  && S' \ar[u]|-{ w} \ar[llu]^-{u'\in \CL} \ar[rru]_-{x'}  
 & &  }  
\end{aligned}
\end{eqnarray}
For $r : V\to U$ in $\CE$, the functor $\overleftarrow{\mathbb{X}}r :\overleftarrow{\mathbb{X}}U\to \overleftarrow{\mathbb{X}}V$ takes the object $(u,x)$ to
$(u_r, (\mathbb{X}r_u)x)$ with notation as in \eqref{cocopbs} which is used
also in defining $\overleftarrow{\mathbb{X}}r$ on morphisms.
The unit $\rho_{\mathbb{X}} : \mathbb{X}\to \overleftarrow{\mathbb{X}}$ for the pseudomonad has components
at $U$ the functor $\rho_U$ taking $x\in \mathbb{X}U$ to $(1_U, x)\in \overleftarrow{\mathbb{X}}U$ and $\zeta : x \to x'$ to $(1_U,\zeta)$. For emphasis we state the dual form of Proposition~\ref{ladjtounit}:
\begin{corollary} Suppose $\CL$ is a protocalibration of the bicategory $\CC^{\mathrm{co}}$. Then
$\mathbb{X}\in \mathrm{CAT}\CC$ is $\CL$-complete if and only if $\rho_{\mathbb{X}}$
has a right adjoint in $\mathrm{CAT}\CC$.   
\end{corollary} 

We think of the right adjoint to $\rho_{\mathbb{X}}$ as assigning a product to the indexed families of generalised objects of $\mathbb{X}$.

\section{$\CR$-Spans and $\CR$-cocompleteness}\label{RSaRc}

We look at slight variants of B\'enabou's bicategory of spans \cite{Ben1967}.
Let $\CR$ be a protocalibration of the bicategory $\CC$; see Definition~\ref{protocalib}. 
A {\em (left) $\CR$-span from $U$ to $V$ in $\CC$} is a diagram $U\xla{u}S\xra{v}V$ in 
$\CC$ with $u\in \CR$. There is a bicategory $\mathrm{Spn}_{\CR}\CC$ whose objects
are those of $\CC$, whose morphisms are the $\CR$-spans, and whose 2-morphisms
(denoted $[f,\sigma]$ by abuse of language) are the isomorphism classes of diagrams \eqref{Rspnmor}.\footnote{It is suggestive to compare with diagram \eqref{Rcocompletion}.
Here the object $\mathbb{X}$ is replaced by a representable $\CC(-,V)$.}
   \begin{eqnarray}\label{Rspnmor}
\begin{aligned}
\xymatrix{
& & S \ar[d]|-{ f} \ar[lld]_-{u\in \CR} \ar[rrd]^-{ v}  
\ar @{} [rd] | {\stackrel{\sigma} \Leftarrow} & &
\\
U    && S' \ar[ll]^-{u'\in \CR} \ar[rr]_-{v'} \ar @{} [lu] | {\stackrel{\cong} \Leftarrow} & & V }  
\end{aligned}
\end{eqnarray} 
The composition of $\CR$-spans is defined on morphisms using bipullback.
To define the composition on 2-morphisms, we use the constructions (in dual form
for groupoid opfibrations rather than fibrations) given in Proposition 5.2 (iii) and
Proposition 8.4 of \cite{134}. 
Explicitly, consider diagram \eqref{componRspnmor}
   \begin{eqnarray}\label{componRspnmor}
\begin{aligned}
\xymatrix{ 
& & S \ar[d]|-{ f} \ar[lld]_-{u\in \CR} \ar[rrd]^-{ v}  
\ar @{} [rd] | {\stackrel{\sigma} \Leftarrow} & &  P  \ar @{} [d] | {\stackrel{\cong} \Leftarrow}\ar[ll]_-{r_v} \ar[rr]^-{v_r}
& & T \ar[d]|-{ g} \ar[lld]_-{r\in \CR} \ar[rrd]^-{ s}  
\ar @{} [rd] | {\stackrel{\sigma} \Leftarrow}
\\
U    && S' \ar[ll]^-{u'\in \CR} \ar[rr]_-{v'} \ar @{} [lu] | {\stackrel{\cong} \Leftarrow} & & V
&& T' \ar[ll]^-{u'\in \CR} \ar[rr]_-{v'} \ar @{} [lu] | {\stackrel{\cong} \Leftarrow} & & W }  
\end{aligned}
\end{eqnarray} 
 where we have bipullbacks \eqref{Lspnbipbs}.
 \begin{eqnarray}\label{Lspnbipbs}
\begin{aligned}
\xymatrix{
P \ar[rr]^-{v_r} \ar[d]_-{r_v} \pullbackcorner && T \ar[d]^-{r} \\
S \ar[rr]_-{v} && V }
\qquad
\xymatrix{
P' \ar[rr]^-{v'_{r'}} \ar[d]_-{r'_{v'}} \pullbackcorner && T' \ar[d]^-{r'} \\
S' \ar[rr]_-{v'} && V}
\end{aligned}
\end{eqnarray}
Using the groupoid opfibration property of $r'$, we obtain a 2-morphism $\kappa : gv_r\Ra k$
and $r'k\cong v'fr_v$ such that $r'\kappa$ is the composite 
$$r'gv_r\cong rv_r\cong vr_v\xra{\sigma r_v} v'fr_v\cong r'k \ .$$
Now define $h : P\to P'$, up to isomorphism using the bipullback property of $P'$, by $r'_{v'}h\cong f r_v$
and $v'_{r'}h\cong k$. Then we have the diagram 
\begin{eqnarray*}
\begin{aligned}
\xymatrix{
& & P \ar[d]|-{ h} \ar[lld]_-{fr_v} \ar[rrd]^-{gv_r}  
\ar @{} [rd] | {\stackrel{\rho} \Leftarrow} & &
\\
S'    && P' \ar[ll]^-{u'\in \CR} \ar[rr]_-{v'_{r'}} \ar @{} [lu] | {\stackrel{\cong} \Leftarrow} & & T' }  
\end{aligned}
\end{eqnarray*} 
 where $\rho$ is the composite $gv_r\xra{\kappa}k\cong v'_{r'}h$.
 Now we define the desired composite of the two span morphisms $[f, \sigma]$ and $[g, \tau]$ 
 appearing in diagram \eqref{componRspnmor} to be $[h, s'\rho]$. 
 
We have a pseudofunctor 
\begin{eqnarray}\label{_*}
(-)_* : \CC\to \mathrm{Spn}_{\CR}\CC
\end{eqnarray}
taking the morphism $U\xra{r} V$ to the $\CR$-span $r_* = (U\xla{1_U}U\xra{r}V)$
and the 2-morphism $\alpha : r\Ra s$ to the 2-morphism
   \begin{eqnarray*}
\begin{aligned}
\xymatrix{
& & U \ar[d]|-{ 1_U} \ar[lld]_-{1_U\in \CR} \ar[rrd]^-{r}  
\ar @{} [rd] | {\stackrel{\alpha} \Leftarrow} & &
\\
U    && U \ar[ll]^-{1_U\in \CR} \ar[rr]_-{s} \ar @{} [lu] | {\stackrel{\cong} \Leftarrow} & & V }  
\end{aligned}
\end{eqnarray*} 
If $r\in \CR$ then $r_*$ has a right adjoint $r^* = (V\xla{r}U\xra{1_U}U)$ in $\mathrm{Spn}_{\CR}\CE$. Every morphism $U\xla{u}S\xra{v}V$ of $\mathrm{Spn}_{\CR}\CE$ decomposes as
$v_*\circ u^*$.    

Restriction along the opposite of \eqref{_*} defines a pseudofunctor
\begin{eqnarray}\label{^*restrict}
\mathrm{CAT} \mathrm{Spn}_{\CR}\CC \to \mathrm{CAT}\CC
\end{eqnarray}
which is pseudomonadic since $(-)_*$ is the identity on objects 
and the restriction has both left and right biadjoints. 
Re-identification of the pseudo-algebras as follows is an objective variant of Lindner's result \cite{Lind}.

\begin{proposition}\label{Rspnclass}
Restriction along the opposite of \eqref{_*} defines a bicategorical equivalence
\begin{eqnarray*}
\mathrm{CAT}(\mathrm{Spn}_{\CR}\CC) \sim \mathrm{CAT}_{\CR}\CC \ .
\end{eqnarray*}
\end{proposition}
\begin{proof}
We shall describe the inverse biequivalence. Each $\mathbb{X}\in \mathrm{CAT}_{\CR}\CC$ extends along the opposite of $(-)_*$ 
to a pseudofunctor $\mathbb{X}_{\mathrm{sp}} : (\mathrm{Spn}_{\CR}\CC)^{\mathrm{op}} \to \mathrm{Cat}$ which
agrees with $\mathbb{X}$ on objects and is defined on morphisms by 
$\mathbb{X}_{\mathrm{sp}}(v_*\circ u^*) = \mathbb{X}_*u\circ \mathbb{X}v$.
For 2-morphisms, we define $\mathbb{X}_{\mathrm{sp}}[f,\sigma] : \mathbb{X}_*u\circ \mathbb{X}v\to \mathbb{X}_*u'\circ \mathbb{X}v'$ to be the mate under the adjunction $\mathbb{X}_*u\dashv \mathbb{X}u$ of the composite
\begin{eqnarray*}
\mathbb{X}v\xra{\mathbb{X}\sigma} \mathbb{X}f\circ\mathbb{X}v'\xra{\mathbb{X}f\circ \eta_{u'}\circ \mathbb{X}v'}\mathbb{X}f\circ \mathbb{X}u'\circ \mathbb{X}_*u'\circ \mathbb{X}v' \cong \mathbb{X}u\circ \mathbb{X}_*u'\circ \mathbb{X}v'
\end{eqnarray*}
where $\eta_{u'}$ is the unit for $\mathbb{X}_*{u'}\dashv \mathbb{X}u'$.
The composition-preservation constraints come from the pointwiseness in $\CR$-cocompleteness: referring to the diagram \eqref{componRspnmor}, we have
\begin{eqnarray*}
\mathbb{X}_{\mathrm{sp}}(v_*\circ u^*)\circ \mathbb{X}_{\mathrm{sp}}(s_*\circ r^*) & \cong & \mathbb{X}_*u\circ \mathbb{X}v\circ \mathbb{X}_*r\circ \mathbb{X}s \\
& \cong & \mathbb{X}_*u\circ \mathbb{X}_*r_v\circ \mathbb{X}v_r\circ \mathbb{X}s \\
& \cong & \mathbb{X}_*(u\circ r_v)\circ \mathbb{X}(s\circ v_r) \\
&  \cong & \mathbb{X}_{\mathrm{sp}}((s\circ v_r)_* \circ (u\circ r_v)^*) \\
&  \cong & \mathbb{X}_{\mathrm{sp}}(s_*\circ {v_r}_* \circ {r_v}^*\circ u^*) \\
& \cong & \mathbb{X}_{\mathrm{sp}}(s_*\circ r^*\circ v_*\circ u^*) \ .  
\end{eqnarray*}
Clearly $\mathbb{X}_{\mathrm{sp}}\circ (-)_*\simeq \mathbb{X}$.

For any $\CR$-cocontinuous pseudonatural transformation $\theta : \mathbb{X}\to \mathbb{Y}$,
we have the diagram
\begin{eqnarray*}
 \begin{aligned}
\xymatrix{
\mathbb{X}V \ar[d]_{\theta_V}^(0.5){\phantom{aaaaa}}="1" 
\ar[rr]^{\mathbb{X}v}  
&& \mathbb{X}S \ar[rr]^{\mathbb{X}_*u} 
     \ar[d]^{\theta_S}_(0.5){\phantom{aaaaa}}="2" \ar@{<=}"1";"2"^-{\theta_v}_-{\cong} 
     \ar@{}[d]^{\phantom{aaaaa}}="3"
&& \mathbb{X}U \ar[d]^{\theta_U}_(0.5){\phantom{aaaaa}}="4" \ar@{<=}"3";"4"^-{\mathring{\theta}_u^{-1}}_-{\cong}
\\
\mathbb{Y}V \ar[rr]_-{\mathbb{Y}v} && \mathbb{Y}S \ar[rr]_-{\mathbb{Y}_*u} && \mathbb{Y}U 
}
 \end{aligned}
\end{eqnarray*}
for all $S\xra{u}U$ in $\CR$ and $S\xra{v}V$ in $\CC$ (using terminology from Proposition~\ref{CB} for the mate of $\theta_u$). The pasted composite diagram exhibits
the extension of $\theta$ to a pseudonatural transformation $\theta_{\mathrm{sp}} : \mathbb{X}_{\mathrm{sp}}\to \mathbb{Y}_{\mathrm{sp}}$. For a modification
$\theta \Ra \phi : \mathbb{X}\to \mathbb{Y}$, we can keep the same components on objects
to define the extended modification $\theta_{\mathrm{sp}} \Ra \phi_{\mathrm{sp}} : \mathbb{X}_{\mathrm{sp}}\to \mathbb{Y}_{\mathrm{sp}}$.   
\end{proof}

Let $\CL$ be a protocalibration of the bicategory $\CC^{\mathrm{co}}$. So we can define 
$$\tensor[_{\CL}]{\mathrm{Spn}}{}\CC = (\mathrm{Spn}_{\CL}\CC^{\mathrm{co}})^{\mathrm{co}} \ .$$
The 2-morphisms are represented by diagrams of the form \eqref{Lspnmor}.
   \begin{eqnarray}\label{Lspnmor}
\begin{aligned}
\xymatrix{
& & S  \ar[lld]_-{u\in \CL} \ar[rrd]^-{ v}  
\ar @{} [rd] | {\stackrel{\kappa} \Rightarrow} & &
\\
U    && S' \ar[u]|-{ f} \ar[ll]^-{u'\in \CL} \ar[rr]_-{v'} \ar @{} [lu] | {\stackrel{\cong} \Rightarrow} & & V }  
\end{aligned}
\end{eqnarray} 

For completeness we state the dual form of Proposition~\ref{Rspnclass}.
\begin{proposition}\label{Lspnclass}
Restriction along the opposite of $(-)_* : \CC \to \tensor[_{\CL}]{\mathrm{Spn}}{}\CC$  defines a bicategorical equivalence
\begin{eqnarray*}
 \mathrm{CAT}(\tensor[_{\CL}]{\mathrm{Spn}}{}\CC) \sim \tensor[_{\CL}]{\mathrm{CAT}}{}\CC \ .
\end{eqnarray*}
\end{proposition}

\section{A distributive law from compatibility}\label{Adlfc}

Let $\CE$ be a category with pullbacks. The functor defined by pulling back along 
$u : S\to U$ is denoted by $\Delta_u : \CE_{/ U}\to \CE_{/ S}$. 
There is always a left adjoint to $\Delta_u$ denoted by $\Sigma_u$ and obtained by
composition with $u$.
Recall that a morphism $u : S\to U$ is called {\em powerful} (or ``exponentiable'') when $\Delta_u$ has a right adjoint $\Pi_u$. Then, for each $a : A\to U$, the counit of the
adjunction $\Delta_u\dashv \Pi_u$ gives a morphism $e$ as in diagram~\eqref{pbpow}. 
\begin{eqnarray}\label{pbpow}
\begin{aligned}
\xymatrix{
A \ar[drr]_-{a}  && P \ar[d]^-{\Delta_u\Pi_ua} \ar[ll]_-{e} \ar[rr]^-{\bar{u}}  \pullbackcorner && B \ar[d]^-{\Pi_ua} \\
  && S \ar[rr]_-{u} && U}
\end{aligned}
\end{eqnarray}
Let $\CP$ denote the set of powerful morphisms in $\CE$;
by Corollary 2.6 of \cite{104}, $\CP$ is a protocalibration in the sense of Definition~\ref{protocalib}.

For this section fix two protocalibrations $\CL$ and $\CR$ of the category $\CE$.
We are particularly interested in the $\CR$-cocompletion and the
$\CL$-completion. 
So, in this section, $\overrightarrow{\mathbb{X}}$ will mean the $\CR$-cocompletion of $\mathbb{X}\in \mathrm{CAT}\CE$ and $\overleftarrow{\mathbb{X}}$ will mean the $\CL$-completion of $\mathbb{X}$. 

Our Definition~\ref{compat} is motivated by concepts in Blumberg-Hill~\cite{BH18, BH22} and Chan~\cite{Ch}.

\begin{definition}\label{compat} Suppose $\CL$ and $\CR$ are protocalibrations of the category $\CE$. We say the pair $(\CL,\CR)$ is {\em compatible} when $\CL\subseteq \CP\cap \CR$ and $\Pi_rv\in \CR$ whenever $r\in \CL$ and $v\in \CR$. 
\end{definition}

\begin{example}
In a category $\CE$ with pullbacks, the maximum example of a compatible pair is $(\CP,\CE)$.
\end{example}

\begin{lemma}\label{Rlowerstar}
Suppose the pair $(\CL,\CR)$ of protocalibrations of $\CE$ is compatible. 
Spans in $\CE$ with left leg invertible and right leg in $\CR$ are morphisms in
$\tensor[_{\CL}]{\mathrm{Spn}}{}\CE$ and as such
form a protocalibration of $\tensor[_{\CL}]{\mathrm{Spn}}{}\CE$ denoted by $\CR_*$.
\end{lemma}
\begin{proof}
$\CR_*$ consists of the spans isomorphic to $w_* = (V\xla{1_V}V\xra{w}V')$ for some 
$V\xra{w}V'$ in $\CR$. The only condition
for a protocalibration which is not totally clear is (G); but this follows from the fact that
each square \eqref{lowerstarbpb} is a bipullback of a groupoid fibration.
\end{proof}
\begin{equation}\label{lowerstarbpb}
\begin{aligned}
\xymatrix{
\tensor[_{\CL}]{\mathrm{Spn}}{}\CE(U,V)^{\mathrm{op}} \ar[rr]^-{\mathrm{incl.}} \ar[d]_-{w_*\circ -} \pullbackcorner
&& \mathrm{Spn}\CE(U,V) \ar[d]^-{w_*\circ -} \\
\tensor[_{\CL}]{\mathrm{Spn}}{}\CE(U,V')^{\mathrm{op}} \ar[rr]_-{\mathrm{incl.}} && \mathrm{Spn}\CE(U,V')}
\end{aligned}
\end{equation}


Suppose $(\CL,\CR)$ is compatible. There is a bicategorical distributive law 
$$\lambda : \overleftarrow{(-)}\circ \overrightarrow{(-)} \lra \overrightarrow{(-)}\circ \overleftarrow{(-)} $$
where the component functor
$$\lambda_U :\overleftarrow{\overrightarrow{\mathbb{X}}}U\lra \overrightarrow{\overleftarrow{\mathbb{X}}}U$$
takes $(U\xla{u}S, S\xla{a}A, x \in \mathbb{X}A)$, where $u\in \CL$ and $a\in \CR$, to 
$$(U\xla{\Pi_ua}B, B\xla{\bar{u}}P, (\mathbb{X}e)x\in \mathbb{X}P)$$ 
in the notation of diagram~\eqref{pbpow}.

\begin{proposition} 
If $\mathbb{X}\in \mathrm{CAT}\CE$ is $\CL$-complete then so is the $\CR$-cocompletion $\overrightarrow{\mathbb{X}}$. Moreover, the unit $\eta_{\mathbb{X}} : \mathbb{X}\to \overrightarrow{\mathbb{X}}$ is in $\tensor[_{\CL}]{\mathrm{CAT}}{}\CE$. 
\end{proposition}
\begin{proof}  
Suppose $V\xra{r}U$ is in $\CL$. We will use the right adjoints $\Pi_r$ and $\mathbb{X}_!r$ available because of powerfulness and of $\CL$-completeness of $\mathbb{X}$.

We must show that $\overrightarrow{\mathbb{X}}r$ has a right adjoint. For each $(V\xla{v}T,y\in \mathbb{X}T)\in \overrightarrow{\mathbb{X}}V$, we must produce an object $(U\xla{v'}R,y'\in \mathbb{X}R)\in \overrightarrow{\mathbb{X}}U$ and a bijection
\begin{eqnarray}\label{!rbij}
\overrightarrow{\mathbb{X}}V((\Delta_ru,(\mathbb{X}\bar{r})x),(v,y)) \cong \overrightarrow{\mathbb{X}}U((u,x),(v',y')) 
\end{eqnarray}
natural in $(u,x)\in \overrightarrow{\mathbb{X}}U$, where the following square is a pullback. 
\begin{equation*}
\xymatrix{
P \ar[rr]^-{\bar{r}} \ar[d]_-{\Delta_ru} \pullbackcorner && S \ar[d]^-{u} \\
V \ar[rr]_-{r} && U}
\end{equation*}
Take an element $(w,\xi)$ as shown in \eqref{eltlhsbij} of the left-hand side of \eqref{!rbij}.
    \begin{eqnarray}\label{eltlhsbij}
\begin{aligned}
\xymatrix{
& & P \ar[d]_{ w} \ar[lld]_-{\Delta_ru} \ar[rrd]^-{(\mathbb{X}\bar{r}) x}  
\ar @{} [rd] | {\stackrel{\xi} \Leftarrow} & &
\\
V  && T \ar[ll]^-{v} \ar[rr]_-{y} & & \mathbb{X} }  
\end{aligned}
\end{eqnarray}
Since $r$ is powerful, morphisms $w$ making the left triangle of \eqref{eltlhsbij} commute are in natural bijection with morphisms $w' : S \to R$ such that $\Pi_rv\circ w' = u$ where $R$ is the domain
of $\Pi_rv$. Since $\bar{r}$ is in $\CL$ and $\mathbb{X}$ is $\CL$-complete, 
morphisms $\xi : (\mathbb{X}\bar{r})x\to (\mathbb{X}w)y$ are in natural bijection with morphisms $\bar{\xi} : x\to (\mathbb{X}_!\bar{r})(\mathbb{X}w)y$. The diagram below shows
how $w$ and $w'$ are related. 
 \begin{eqnarray*}
 \begin{aligned}
\xymatrix{ && P \ar[dll]_-{w} \ar[rr]^-{\bar{r}} \ar[d]^-{\Delta_{\tilde{r}}w'} \pullbackcorner && S \ar[d]^-{w'}    \\
T \ar[drr]_-{v}  && Q \ar[d]^-{\Delta_r\Pi_rv} \ar[ll]_-{e} \ar[rr]^-{\tilde{r}}  \pullbackcorner && R \ar[d]^-{\Pi_rv} \\
  && V \ar[rr]_-{r} && U}
\end{aligned}
\end{eqnarray*}
Also, we have the Chevalley-Beck
isomorphism $(\mathbb{X}w') (\mathbb{X}_!\tilde{r}) \cong (\mathbb{X}_!\bar{r})(\mathbb{X}\Delta_{\tilde{r}}w')$. Consequently, we have
$(\mathbb{X}_!\bar{r})(\mathbb{X}w) \cong (\mathbb{X}_!\bar{r})(\mathbb{X}\Delta_{\tilde{r}}w')(\mathbb{X}e)\cong (\mathbb{X}w') (\mathbb{X}_!\tilde{r})(\mathbb{X}e)$. So
$\bar{\xi} : x\to (\mathbb{X}_!\bar{r})(\mathbb{X}w)y$ are in bijection with
$\bar{\xi} : x\to (\mathbb{X}w') (\mathbb{X}_!\tilde{r})(\mathbb{X}e)y$.
Therefore we have the desired bijection \eqref{!rbij} with $v'=\Pi_rv$ and $y'= (\mathbb{X}_!\tilde{r})(\mathbb{X}e)y$. 

Since $\overrightarrow{\mathbb{X}}$ is $\CR$-cocomplete, it already satisfies the CB-condition for
left adjoints $\overrightarrow{\mathbb{X}}_*r$. By replacing all functors in that invertibility 
condition by their right adjoints when $r\in \CL$, we induce an invertible mate yielding
the CB-condition for the right adjoints $\overrightarrow{\mathbb{X}}_!r$. 
So $\overrightarrow{\mathbb{X}}$ is $\CL$-complete.   

The second sentence of the proposition follows on checking that the canonical natural
transformation
\begin{eqnarray*}
\begin{aligned}
  \xymatrix{
\mathbb{X}V \ar[d]_{\mathbb{X}_!r}^(0.5){\phantom{aaaaa}}="1" \ar[rr]^{\eta_{\mathbb{X}}V}  && \overrightarrow{\mathbb{X}}V \ar[d]^{\overrightarrow{\mathbb{X}}_!r}_(0.5){\phantom{aaaaa}}="2" \ar@{=>}"1";"2"^-{ }
\\
\mathbb{X}U \ar[rr]_-{\eta_{\mathbb{X}}U} && \overrightarrow{\mathbb{X}}U 
}
 \end{aligned}
\end{eqnarray*}
is invertible. 
\end{proof}

The pseudomonad $\overrightarrow{(-)}$ on $\mathrm{CAT}\CE$ therefore lifts to $\tensor[_{\CL}]{\mathrm{CAT}}{}\CE$. 
\begin{eqnarray}\label{lifttoL}
\begin{aligned}
  \xymatrix{
\tensor[_{\CL}]{\mathrm{CAT}}{}\CE \ar[d]_{\mathrm{incl}}^(0.5){\phantom{aaaaaaaa}}="1" \ar[rr]^{\tensor[_{\CL}]{\overrightarrow{(-)}}{}}  && \tensor[_{\CL}]{\mathrm{CAT}}{}\CE \ar[d]^{\mathrm{incl}}_(0.5){\phantom{aaaaaaaa}}="2" \ar@{=>}"1";"2"^-{\simeq}
\\
\mathrm{CAT}\CE \ar[rr]_-{\overrightarrow{(-)}} && \mathrm{CAT}\CE 
}
 \end{aligned}
\end{eqnarray}
The pseudo-algebras for the lifted pseudomonad are those of the composite pseudomonad
$\mathbb{X}\mapsto \overrightarrow{\overleftarrow{\mathbb{X}}}$ on $\mathrm{CAT}\CE$ 
form the 2-category we denote by $\tensor[_{\CL}]{\mathrm{CAT}}{_{\CR}}\CE$ at the top
of the diamond \eqref{adjointsquare} of inclusion 2-functors.
Using Lemma~\ref{Rlowerstar} and Proposition~\ref{Rspnclass}, we have biequivalences   
\begin{eqnarray*}
\tensor[_{\CL}]{\mathrm{CAT}}{_{\CR}}\CE \sim \mathrm{CAT}_{\CR_*}(\tensor[_{\CL}]{\mathrm{Spn}}{}\CE)\sim \mathrm{CAT}(\mathrm{Spn}_{\CR_*}(\tensor[_{\CL}]{\mathrm{Spn}}{}\CE)) \ . 
\end{eqnarray*}
While the composite pseudomonad
is not Kock-Z\"oberlein or dual Kock-Z\"oberlein, the pseudo-algebra structures are unique up to isomorphism: it is a property of an $\CE$-variable category to be a pseudo-algebra for both $\overleftarrow{(-)}$ and $\overrightarrow{(-)}$, and to satisfy the distributive law.   
\begin{eqnarray}\label{adjointsquare}
\begin{aligned}
\xymatrix{
& \tensor[_{\CL}]{\mathrm{CAT}}{_{\CR}}\CE \ar[rd]^-{ } \ar[ld]^-{}  & \\
\tensor[_{\CL}]{\mathrm{CAT}}{}\CE   \ar[rd]_-{} & & \mathrm{CAT}_\CR\CE \ar[ld]^-{ } \\
& \mathrm{CAT}\CE }
 \end{aligned}
\end{eqnarray}

\section{Polynomials as spans of spans}\label{Pasos} 
 Let $\CE$ be a category with pullbacks. 

There is a bicategory $\mathrm{Ply}\CE$ of polynomials in $\CE$; in the notation of 
Walker \cite{Walker2019} it is $\mathrm{Poly}(\CE)$ and in the notation of Proposition 8.4 of \cite{134} it is the bicategory $\mathrm{PolySpn}{\CE}$.
The objects are those of $\CE$. 
The morphisms from $X$ to $Y$ are {\em polynomials} $p = (X\xleftarrow{r} A\xrightarrow{n}B\xrightarrow{t} Y)$ in $\CE$
with $n$ powerful. 
The morphisms between polynomials are as in Subsection 2.11 of \cite{GambKock} and Definition 8.2 of \cite{134}: they are isomorphism classes of commutative diagrams
\begin{equation}\label{morph_of_polys}
 \begin{aligned}
\xymatrix{
& A \ar[dl]_-{r}   \ar[r]^-{n}  &   B  \ar[dr]^{t}  &  \\
X & P \ar[d]|-{g_{n'}} \ar[u]_-{\ell}  \ar[r]|-{n'_g} \pullbackcorner & B \ar[d]_-{g}  \ar[u]^-{1_B} & Y   \\
& A' \ar[ul]^-{r'} \ar[r]_-{n'} & B'   \ar[ur]_-{t'}  & 
}
 \end{aligned}
\end{equation}
wherein the square as indicated is a pullback. Composition uses bipullback of spans of spans (see \cite{GambKock, Weber2015, 134}). 
As a consequence of Theorem~\ref{PlyasSpnSpn} and Proposition~\ref{fincoprodsglobloc} below, finite coproducts exist in each $\mathrm{Ply}{\CE}(X,Y)$ and are preserved by precomposition with morphisms, yet generally not by postcomposition  
making the bicategory $\mathrm{Ply}{\CE}$ 
enriched in a certain skew monoidal bicategory $\mathrm{Cat}_{+ \mathrm{sk}}$
(extending a bit the terminology of \cite{GarLem}).

Recall the philosophy of \cite{134} that it is useful to regard $p$ as a span
\begin{eqnarray}\label{spn_of_spns}
p = (X\xleftarrow{(n,A,r)} B\xrightarrow{(1_B,B,t)} Y)
\end{eqnarray}
in $\mathrm{Spn}\CE$.  Then the morphisms interpret as diagrams
  \begin{eqnarray}\label{strongpolymorph}
\begin{aligned}
\xymatrix{
& & B \ar[d]|-{ (1, B, g)} \ar[lld]_-{(n,A,r)} \ar[rrd]^-{ (1,B,t)} \ar @{} [ld] | {\stackrel{\lambda} \Leftarrow} 
\ar @{} [rd] | {\stackrel{\phantom{a} } =} & &
\\
X   && B' \ar[ll]^-{(n',A',r')} \ar[rr]_-{(1,B',t')} & & Y  \ .}
\end{aligned}
\end{eqnarray}

We modify this a bit by taking a compatible pair $(\CL,\CR)$ of protocalibrations on $\CE$
and restricting the polynomials \eqref{spn_of_spns} to those with $n\in \CL$ and $t\in \CR$.
Restricting the morphisms of $\mathrm{Ply}{\CE}$ in this way, we obtain a sub-bicategory
which we denote by $\tensor[_{\CL}]{\mathrm{Ply}}{_{\CR}}\CE$. 

\begin{theorem}\label{PlyasSpnSpn} If $(\CL,\CR)$ is a compatible pair of protocalibrations on the category $\CE$ then
$$(\tensor[_{\CL}]{\mathrm{Ply}}{_{\CR}}\CE)^{\mathrm{op}} = \mathrm{Spn}_{\mathrm{\CR_*}}(\tensor[_{\CL}]{\mathrm{Spn}}{}\CE) \ .$$
\end{theorem}
\begin{proof}
The 2-morphisms of $\mathrm{Spn}_{\CR_*}(\tensor[_{\CL}]{\mathrm{Spn}}{}\CE)$ 
are isomorphism classes of morphisms as pictured in \eqref{pspolymorph}.  
 \begin{eqnarray}\label{pspolymorph}
\begin{aligned}
\xymatrix{
& & B \ar[d]|-{ (1, B, g)} \ar[rrd]^-{(n,A,r)} \ar[lld]_-{ (1,B,t)} \ar @{} [ld] | {\stackrel{\cong } \Leftarrow} 
\ar @{} [rd] | {\stackrel{\lambda} \Leftarrow}  & &
\\
Y   && B' \ar[rr]_-{(n',A',r')} \ar[ll]^-{(1,B',t')} & & X  \ .}
\end{aligned}
\end{eqnarray}
Unpacking this shows that the data involved is exactly as in
the diagram \eqref{morph_of_polys} in $\CE$.
While unpacking, just note that the 2-morphism $\lambda$ in
$\tensor[_{\CL}]{\mathrm{Spn}}{}\CE$ is in the reverse direction as a morphism of spans. 
\end{proof}
\begin{corollary}\label{Walker} $\mathrm{Ps}(\tensor[_{\CL}]{\mathrm{Ply}}{_{\CR}}\CE, \mathrm{Cat}) \sim \tensor[_{\CL}]{\mathrm{CAT}}{_{\CR}}\CE$
\end{corollary}
\begin{proof} Using Theorem~\ref{PlyasSpnSpn} and Propositions~\ref{Rspnclass} and \ref{Lspnclass}, we have
\begin{eqnarray*}
\phantom{aaaaaaaaaa} \mathrm{Ps}(\tensor[_{\CL}]{\mathrm{Ply}}{_{\CR}}\CE, \mathrm{Cat}) &\sim & \mathrm{CAT}(\mathrm{Spn}_{\mathrm{\CR_*}}(\tensor[_{\CL}]{\mathrm{Spn}}{}\CE)) \\
&\sim & \mathrm{CAT}_{\CR_*}(\tensor[_{\CL}]{\mathrm{Spn}}{}\CE) \\
& \sim & \tensor[_{\CL}]{\mathrm{CAT}}{_{\CR}}\CE \ . \ \phantom{aaaaaaaaaaaaaaaa} \qedhere
\end{eqnarray*} 
\end{proof}
\begin{eqnarray}\label{inducingdiamond}
\begin{aligned}
\xymatrix{
& \CE \ar[rd]^-{ (-)_*} \ar[ld]_-{(-)_*}  & \\
\tensor[_{\CL}]{\mathrm{Spn}}{}\CE   \ar[rd]_-{ } & & \mathrm{Spn}_\CR\CE \ar[ld]^-{ } \\
& (\tensor[_{\CL}]{\mathrm{Ply}}{_{\CR}}\CE)^{\mathrm{op}} }
 \end{aligned}
\end{eqnarray}

\section{Finite coproducts for spans and polynomials}\label{Fcfsap}
Examples of the concepts of this section are provided by lextensive categories.
Recall that a category $\CE$ is {\em lextensive} (in the sense used by Schanuel 
\cite{Schanuel}) when it has finite limits, finite coproducts, and, equivalences 
$$\mathbf{1} \simeq \CE_{/0} \ \text{ and } \ \CE_{/U}\times \CE_{/V}\simeq \CE_{/U+V} \ .$$ 
The second equivalence takes an object $(A\xra{h}U, B\xra{k}V)\in \CE_{/U}\times \CE_{/V}$
to the object $(A+B\xra{h+k}U+V)\in \CE_{/U+V}$; the fact that this inverse equivalence is
fully faithful means that a commutative triangle
\begin{eqnarray}\label{lextfact}
\begin{aligned}
\xymatrix{
A+B \ar[rd]_{h+k}\ar[rr]^{f}   && A'+B' \ar[ld]^{h'+k'} \\
& U+V  &
}
\end{aligned}
\end{eqnarray}
implies there are unique $r : A\to A'$ and $s : B\to B'$ such that $f = r+s$, $h = h'r$
and $h'=k's$. 

It follows that the initial object $0$ is strictly initial and, for every coproduct diagram $U\xra{i}U+V\xla{j}V$ in $\CE$ and every $R\xra{f}U+V$, there is a diagram
    \begin{eqnarray}\label{lextpbs}
\begin{aligned}
\xymatrix{
S \ar[rr]^-{\bar{i}} \ar[d]_-{h}  && R  \ar[d]_-{f}  && T \ar[d]^-{k} \ar[ll]_-{\bar{j}} \\
U \ar[rr]_-{i} && U+V   && V \ar[ll]^-{j}}
\end{aligned}
\end{eqnarray}
in which the squares are pullbacks and the top row is also a coproduct diagram.

This should motivate the following bicategorical concepts.
 
Suppose the bicategory $\CC$ admits the binary bicoproduct 
$U\xra{i}U+V\xla{j}V$ and $\CR$ is a protocalibration of $\CC$. 
For the pseudofunctor $\mathbb{C}_{\CR} : \CC^{\mathrm{op}}\to \mathrm{Cat}$ of 
Definition~\ref{C_R}, we have a canonical comparison functor
$$\mathbb{C}_{\CR}(U+V)\to \mathbb{C}_{\CR}U\times \mathbb{C}_{\CR}V$$ 
taking $Z\xra{w}U+V$ to the pair $(X\xra{w_i}U, Y\xra{w_j}V)$ defined by bipullbacks
\begin{eqnarray*}
\xymatrix{
X\pullbackcorner \ar[d]_{w_i}^(0.5){\phantom{aaaaaa}}="1" \ar[rr]^{i_w}  && Z \ar[d]^{w}_(0.5){\phantom{aaaaaa}}="2" \ar@{=>}"1";"2"^-{\cong}
\\
U \ar[rr]_-{i} && U+V 
}
\qquad
\xymatrix{
Y\pullbackcorner \ar[d]_{w_j}^(0.5){\phantom{aaaaaa}}="1" \ar[rr]^{j_w}  && Z \ar[d]^{w}_(0.5){\phantom{aaaaaa}}="2" \ar@{=>}"1";"2"^-{\cong}
\\
V \ar[rr]_-{j} && U+V \ . 
}
\end{eqnarray*}

\begin{lemma}\label{adjforcancomp} 
Suppose the bicategory $\CC$ admits, and the protocalibration $\CR$ is closed under, finite bicoproducts. Then the canonical comparison functors
\begin{eqnarray*}
\mathbb{C}_{\CR}0 \to \mathbf{1} \ \text{  and  } \ \mathbb{C}_{\CR}(U+V)\to \mathbb{C}_{\CR}U\times \mathbb{C}_{\CR}V
\end{eqnarray*}
have left adjoints respectively defined on objects by
$$ (0\in \mathbf{1}) \mapsto (0\xra{1_0}0) \ \text{ and } \ (X\xra{u\in \CR}U, Y\xra{v\in \CR}V)\mapsto (X+Y\xra{u+v\in \CR}U+V) \ .$$ 
\end{lemma}
\begin{proof}
This is an exercise in using the universal properties of bicoproduct and bipullback.
\end{proof}
\begin{definition} Suppose $\CR$ is a protocalibration of the bicategory $\CC$.
We say a bicategory $\CC$ is {\em $\CR$-extensive} when it admits finite bicoproducts and these are taken into finite 
biproducts by the pseudofunctor $\mathbb{C}_{\CR} : \CC^{\mathrm{op}}\to \mathrm{Cat}$
of Definition~\ref{C_R}.
\begin{eqnarray*}
\mathbb{C}_{\CR}0 \simeq \mathbf{1} \ \text{  and  } \ \mathbb{C}_{\CR}(U+V)\simeq \mathbb{C}_{\CR}U\times \mathbb{C}_{\CR}V
\end{eqnarray*}
When $\CR$ is closed under finite bicoproducts, the inverse equivalences are 
the left adjoints in Definition~\ref{adjforcancomp}. 
\end{definition}
\begin{example}\label{lextensiveE=R} Suppose $\CE$ is a lextensive category. Then the maximum protocalibration $\CE$ is closed under finite coproducts and $\CE$ is $\CE$-extensive. 
\end{example}

\begin{example} Again suppose $\CE$ is a lextensive category. The subcategory $\CP$ consisting of the powerful morphisms is closed under finite coproducts. To see this, first realize that injections $X\xra{i}X+Y$ are powerful since $\CE_{/ (X+Y)}\xra{\Delta_i}\CE_{/ X}$
transports across equivalence to the projection $\CE_{/ X}\times \CE_{/ Y}\to \CE_{/ X}$ which has a right adjoint since $\CE_{/ Y}$ has a terminal object. Then notice that pullback along $0\to A$ is equivalent to the functor $\CE_{/A}
\to \mathbf{1}$ whose right adjoint picks out the terminal object and, if $A\xra{u}X$ and $B\xra{v}X$ are powerful, then pullback along $A+B\xra{[u,v]} X$ is equivalent to
$\CE_{/X}\xra{(\Delta_u,\Delta_v)}\CE_{/A}\times \CE_{/B}$ which has a right
adjoint taking $(U\xra{a}A, V\xra{b}B)$ to the product $\Pi_u(a)\times\Pi_v(b)$ in $\CE_{/X}$. 
Indeed, $\CE$ is $\CP$-extensive because the equivalence 
$\CE_{/ (X+Y)}\simeq \CE_{/ X}\times \CE_{/ Y}$ restricts to the powerful objects.   
\end{example}

\begin{proposition}\label{Rcocomrestricts} Suppose the bicategory $\CC$ is $\CR$-extensive and that the protocalibration $\CR$ is closed under the finite bicoproducts.
If $\mathbb{X}\in \mathrm{CAT}\CC$ preserves finite bicategorical products then the $\CR$-cocompletion $\overrightarrow{\mathbb{X}}$ also preserves finite bicategorical products.  
\end{proposition}
\begin{proof}
The objects of $\overrightarrow{\mathbb{X}}(U+V)$ are pairs 
$(U+V\xla{w}R, z\in \mathbb{X}R)$ with $w\in \CR$.
Since $\CC$ is $\CR$-extensive, $R\xra{w}U+V$ is equivalent over $U+V$ to some
$S+T\xra{u+v}U+V$ and the morphisms $u$ and $v$ are in $\CR$.
Since $\mathbb{X}$ preserves finite biproducts, $z\in \mathbb{X}R$
gives $(x,y)\in \mathbb{X}S\times \mathbb{X}T$. The assignment $(R,z)\mapsto ((S,x), (T,y))$ is the required equivalence 
$\overrightarrow{\mathbb{X}}(U+V)\simeq \overrightarrow{\mathbb{X}}U\times\overrightarrow{\mathbb{X}}V$. It is clear that $\overrightarrow{\mathbb{X}}0\simeq \mathbf{1}$ since all $0\xla{u\in \CR}S$ are isomorphic equivalences and $\mathbb{X}0\simeq \mathbf{1}$.    
\end{proof}

\begin{proposition}\label{fincoprodsglobloc} Suppose the bicategory $\CC$ is $\CR$-extensive and that the protocalibration $\CR$ is closed under the finite bicoproducts. Then 
\begin{itemize}
\item[i.] the bicategory $\mathrm{Spn}_{\CR}\CC$ admits finite bicoproducts, and 
\item[ii.] each homcategory $\mathrm{Spn}_{\CR}\CC(X,Y)$ admits finite
coproducts preserved by composition on the diagrammatic right by any morphism $Y\xra{(u,S,v)} Z$ where $u\in \CR$. 
\end{itemize} 
\end{proposition}
\begin{proof} The bi-initial object of $\mathrm{Spn}_{\CR}\CC$ is that of $\CC$ denoted $0$. For the 1-morphism part of the universal property: if $0\xla{u}S\xra{v}Y$ is a left $\CR$-span then $u$ is an equivalence (using $\mathbb{C}_{\CR}0 \sim \mathbf{1}$); so $S$ is bi-initial and we conclude that all such spans are isomorphic in
$\mathrm{Spn}_{\CR}\CC$ while the 2-morphism property of $0$ follows from that of $0\in \CC$.

We will show that a bicoproduct $U\xra{i}U+V\xla{j}V$ in $\CC$ gives one $U\xra{i_*}U+V\xla{j_*}V$ in $\mathrm{Spn}_{\CR}\CC$. Note that coprojections $i$, $j$ and codiagonals are in $\CR$. Take spans $U\xla{u}S\xra{v}W$ and $V\xla{r}T\xra{v}W$ with
$u, r\in \CR$ (and hence $u+r \in \CR$). We obtain $U+V\xla{u+r\in \CR}S+T\xra{[v,s]}W$ whose composites with $i_*$ and $j_*$
are isomorphic to the given two spans because our assumptions give the bipullbacks
\begin{eqnarray*}
\xymatrix{
S\pullbackcorner \ar[d]_{u}^(0.5){\phantom{aaaaaa}}="1" \ar[rr]^{\bar{i}}  && S+T \ar[d]^{u+r}_(0.5){\phantom{aaaaaa}}="2" \ar@{=>}"1";"2"^-{\cong}
\\
U \ar[rr]_-{i} && U+V 
}
\qquad
\xymatrix{
T\pullbackcorner \ar[d]_{r}^(0.5){\phantom{aaaaaa}}="1" \ar[rr]^{\bar{j}}  && S+T \ar[d]^{u+r}_(0.5){\phantom{aaaaaa}}="2" \ar@{=>}"1";"2"^-{\cong}
\\
V \ar[rr]_-{j} && U+V \ . 
}
\end{eqnarray*}
The 2-morphism part of the universal property uses
$\CR$-extensivity and the properties of $\CR$:
   \begin{eqnarray*}
\begin{aligned}
\xymatrix{
& & S+T \ar[d]|-{ f+g} \ar[lld]_-{u+r\in \CR} \ar[rrd]^-{ [v,s]}  
\ar @{} [rd] | {\stackrel{[\sigma,\tau]} \Leftarrow} & &
\\
U+V    && S'+T' \ar[ll]^-{u'+r'\in \CR} \ar[rr]_-{[v',s']} \ar @{} [lu] | {\stackrel{\cong} \Leftarrow} & & W \ . }  
\end{aligned}
\end{eqnarray*} 

The coproduct of $X\xla{u\in\CR}S\xra{v}Y$ and $X\xla{u'\in\CR}S'\xra{v'}Y$ in $\mathrm{Spn}_{\CR}\CC(X,Y)$
is $X\xla{[u,u']\in\CR}S+S'\xra{[v,v']}Y$. Contemplation of the diagram 
\begin{eqnarray*}
 \begin{aligned}
\xymatrix{
& & \ar @{} [dd] | {\stackrel{\cong } \Longleftarrow} R+R' \ar[ld]_-{a+b\in \CR} \ar[rd]^-{[c,d]} & & \\
& S+S' \ar[ld]_-{[u,u']\in \CR} \ar[rd]_-{[v,v']} & & T \ar[ld]^-{r\in \CR} \ar[rd]^-{s} & \\
X & & Y & & Z}
 \end{aligned}
\end{eqnarray*}
(where the diamond is a bipullback) reveals the claimed right preservation property.          
\end{proof}

For any bicategory $\CD$ admitting finite bicoproducts, we write $\mathrm{fpCAT}\CD$
for the full sub-bicategory of $\mathrm{CAT}\CD$ consisting of the finite biproduct
preserving pseudofunctors $\CD^{\mathrm{op}}\to \mathrm{Cat}$. If $\CT$ is a
protocalibration of $\CD$, we write $\mathrm{fpCAT}_{\CT}\CD$ for the intersection
of $\mathrm{fpCAT}\CD$ and $\mathrm{CAT}_{\CT}\CD$.    

\begin{proposition}\label{RspnclassC}
Under the hypotheses of Proposition~\ref{fincoprodsglobloc}, the biequivalence of Proposition~\ref{Rspnclass} restricts to a biequivalence
\begin{eqnarray*}
\mathrm{fpCAT}(\mathrm{Spn}_{\CR}\CC) \sim \mathrm{fpCAT}_{\CR}\CC \ .
\end{eqnarray*} 
\end{proposition}
\begin{proof} The pseudofunctor $(-)_* : \CC \to \mathrm{Spn}_{\CR}\CC$ preserves
finite bicategorical coproduct so restriction along its opposite restricts to a pseudofunctor 
from left to right which is an equivalence on homcategories. 
Also, if $\mathbb{X} \in \mathrm{CAT}_{\CR}\CC$
preserves finite bicategorical products then so does its extension $\mathbb{X}_{\mathrm{sp}}$
to $(\mathrm{Spn}_{\CR}\CC)^{\mathrm{op}}$ since $\CR$ is closed under finite bicoproduct.   
\end{proof}

\begin{example}\label{lextensiveE=Rcont} Continuing with Example~\ref{lextensiveE=R}, Proposition~\ref{fincoprodsglobloc} tells us that $\mathrm{Spn}\CE = \mathrm{Spn}_{\CE}\CE$ 
has bicoproducts and local coproducts preserved by right composition. However, 
$\mathrm{Spn}\CE = (\mathrm{Spn}\CE)^{\mathrm{op}}$, so local coproducts are
preserved by composition on the left as well. Furthermore, if $\CR$ is any protocalibration
of the lextensive category $\CE$ closed under finite coproducts then $\CE$ is $\CR$-extensive.
The bicoproducts and local coproducts in $\mathrm{Spn}_{\CR}\CE$ are preserved by the inclusion pseudofunctor into $\mathrm{Spn}\CE$. It follows that the local coproducts in
$\mathrm{Spn}_{\CR}\CE$ are preserved by composition on both sides. 
Consider the symmetric closed monoidal bicategory $\mathrm{Cat}_{+}$ of categories with finite coproducts and 
finite-coproduct-preserving functors; the internal hom $[\CA,\CB]_{+}$ of 
$\mathrm{Cat}_{+}$ is $\mathrm{Cat}_{+}(\CA,\CB)$ with pointwise coproduct
and the unit object is the category $\mathrm{set}$ of finite sets.
Cartesian product $A\times B$ of categories is a direct sum in 
$\mathrm{Cat}_{+}$ and the terminal category $\mathbf{1}$ is also initial. 
What we have is that $\mathrm{Spn}\CE$ is $\mathrm{Cat}_{+}$-enriched and
to say a pseudofunctor $(\mathrm{Spn}\CE)^{\mathrm{op}}\to\mathrm{Cat}_{+}$ is
biproduct preserving is to say it is $\mathrm{Cat}_{+}$-enriched.     
\end{example}

\begin{lemma}\label{Rextensivity} 
Suppose $(\CL,\CR)$ is a compatible pair of protocalibrations on an $\CR$-extensive category $\CE$. Suppose both $\CL$ and $\CR$ are closed under finite coproducts in $\CE$.
The protocalibration $\CR_*$ of 
the bicategory $\tensor[_{\CL}]{\mathrm{Spn}}{}\CE$ (see Lemma~\ref{Rlowerstar}) is 
closed under bicoproducts and $\tensor[_{\CL}]{\mathrm{Spn}}{}\CE$ is $\CR_*$-extensive. 
\end{lemma}
\begin{proof} For this proof, put $\CC = \tensor[_{\CL}]{\mathrm{Spn}}{}\CE$. The injections into
a bicoproduct and the codiagonals are in $\CR_*$. The bicoproduct of $f_*$ and $f'_*$ is isomorphic to $(f+f')_*$ and so in $\CR_*$. It remains to prove $\mathbb{C}_{\CR_*}$ preserves biproducts. 
Up to equivalence, we can take the objects of $\mathbb{C}_{\CR_*}U$ to be 
morphisms $X\xra{a}U$ in $\CR$ since every morphism $X\to U$ in $\CR_*$
is isomorphic to some $X\xra{a_*}U$ with $a\in \CR$. Moreover, every morphism
\begin{eqnarray*}
\xymatrix{
X \ar[rd]_{a_*}^(0.5){\phantom{a}}="1" \ar[rr]^{[u,S,v]}  && Y \ar[ld]^{b_*}_(0.5){\phantom{a}}="2" \ar@{=>}"1";"2"^-{\cong}
\\
& U 
}
\end{eqnarray*}
in $\mathbb{C}_{\CR_*}U$ has a morphism of the form $f_*$ in the isomorphism class 
$[u,S,v]$. Thus we have an equivalence of categories 
$\mathbb{C}_{\CR_*}U\simeq \mathbb{E}_{\CR}U$.\footnote{recalling Definition~\ref{C_R} although now we have $\mathbb{E}$ instead of $\mathbb{C}$ since we have $\CE$ instead of $\CC$.}
Also, the canonical functor 
$\mathbb{C}_{\CR_*}(U+V)\to \mathbb{C}_{\CR_*}U \times \mathbb{C}_{\CR_*}V$ transports across the equivalences to the lextensivity equivalence
$\mathbb{E}_{\CR}(U+V) \simeq \mathbb{E}_{\CR}U\times \mathbb{E}_{\CR}V$.    
\end{proof}

As a corollary of Theorem~\ref{PlyasSpnSpn}, Proposition~\ref{fincoprodsglobloc} and
Lemma~\ref{Rextensivity}, we obtain:

\begin{theorem}\label{plycoprodthm} Under the hypotheses of Lemma~\ref{Rextensivity}, the bicategory $\tensor[_{\CL}]{\mathrm{Ply}}{_\CR}\CE$
admits finite biproducts and its homs admit finite coproducts preserved by composition on the diagrammatic left.
\end{theorem}

\section{Objective Mackey functors}   
 
{\em Suppose the bicategory $\CC$ is finitely bicocomplete and $\CR$-extensive for a protocalibration $\CR$ which, when regarded as a wide subbicategory of $\CC$, is closed under finite bicoproducts.}
 
\begin{definition}
An {\em objective $\CR$-Mackey functor} on $\CC$ is an $\CR$-cocomplete pseudofunctor 
$\mathbb{M} : \CC^{\mathrm{op}} \to \mathrm{Cat}$ which preserves finite biproducts (in the bicategorical sense). 
We define the bicategory of objective $\CR$-Mackey functors by $\mathrm{OMky}_{\CR}\CC = \mathrm{fpCAT}_{\CR}\CC$. By Proposition~\ref{RspnclassC}, 
\begin{eqnarray}\label{RspnclassCM}
\mathrm{fpCAT}(\mathrm{Spn}_{\CR}\CC) \sim \mathrm{OMky}_{\CR}\CC \ .
\end{eqnarray} 
When $\CR = \CC$, we drop the subscripts $\CC$ and refer to {\em objective Mackey functors}.
\end{definition}
\begin{example}\label{oRBf} 
In the terminology of Definition~\ref{C_R}, the object
$\mathbb{C}_{\CR} \in \mathrm{CAT}_{\CR}\CC$ is an objective $\CR$-Mackey functor and, as suggested by the case $\CC = \mathrm{set}^G$ of finite sets acted on by a finite group $G$, is called the
{\em objective $\CR$-Burnside functor}.  
\end{example}
\begin{example}\label{inthomMky}In the terminology of Proposition~\ref{inthom}, if $\mathbb{Z}$ is an objective $\CR$-Mackey functor, so is $[\mathbb{Y},\mathbb{Z}]_{\CR}$. For, we have 
\begin{eqnarray*}
[\mathbb{Y},\mathbb{Z}]_{\CR}(W+W') & \simeq & \mathrm{CAT}_{\CR}\CE(\mathbb{Y},\mathbb{Z}(-\times (W+W')) \\
& \simeq & \mathrm{CAT}_{\CR}\CE(\mathbb{Y},\mathbb{Z}(-\times W)\times \mathbb{Z}(-\times W')) \\
& \simeq & \mathrm{CAT}_{\CR}\CE(\mathbb{Y},\mathbb{Z}(-\times W)) \times \mathrm{CAT}_{\CR}\CE(\mathbb{Y},  \mathbb{Z}(-\times W')) \\
& \simeq & [\mathbb{Y},\mathbb{Z}]_{\CR}W\times [\mathbb{Y},\mathbb{Z}]_{\CR}W' \ .
\end{eqnarray*}
(We are using the last sentence of Corollary~\ref{pseudomonadicity} in the third step.) 
\end{example}
As in Example~\ref{lextensiveE=Rcont}, we consider the symmetric closed monoidal bicategory $\mathrm{Cat}_{+}$ of categories with finite coproducts and 
finite-coproduct-preserving functors; the internal hom $[\CA,\CB]_{+}$ of 
$\mathrm{Cat}_{+}$ is $\mathrm{Cat}_{+}(\CA,\CB)$ with pointwise coproduct.
Cartesian product $A\times B$ of categories is a direct sum in $\mathrm{Cat}_{+}$ and the terminal category $\mathbf{1}$ is also initial. 
\begin{proposition}\label{Catsub+fact}
Every objective $\CR$-Mackey functor $\mathbb{M} : \CC^{\mathrm{op}} \to \mathrm{Cat}$ factors through the forgetful 2-functor $\mathrm{Cat}_{+}\to \mathrm{Cat}$ by a pseudofunctor also denoted by $\mathbb{M}$.  
\end{proposition}
\begin{proof}
By assumption, the codiagonal $\nabla : U+U\to U$ and unique $! : 0\to U$ are in $\CR$. The diagonal of $\mathbb{M}U$ is isomorphic to the composite 
$$\mathbb{M}U \xra{\mathbb{M}\nabla} \mathbb{M}(U+U)\simeq  \mathbb{M}U\times \mathbb{M}U$$
and so has a left adjoint by Proposition~\ref{CB}. Thus $\mathbb{M}U$ has binary coproducts.
Similarly, $\mathbb{M}U \xra{\mathbb{M}!} \mathbb{M}0\simeq \mathbf{1}$ has a left adjoint so $\mathbb{M}U$ has an
initial object. For $f : U\to V$, the functor $\mathbb{M}f$ preserves these finite coproducts as can be seen by applying the CB property to the pullbacks
\begin{eqnarray*}
\xymatrix{
U+U \ar[rr]^-{\nabla} \ar[d]_-{f+f} && U \ar[d]^-{f} \\
V+V \ar[rr]_-{\nabla} && V}
\quad
\xymatrix{ \\ 
& , & }
\quad
\xymatrix{
0 \ar[rr]^-{!} \ar[d]_-{1_0} && U \ar[d]^-{f} \\
0 \ar[rr]_-{!} && V}
\end{eqnarray*}
in our $\CR$-extensive bicategory $\CC$.     
\end{proof}
\begin{corollary}
If $\mathbb{M}$ is an objective $\CR$-Mackey functor and $U\xra{i} U+V\xla{j}V$ is
a coproduct in $\CC$ then
\begin{eqnarray*}
\xymatrix @R-3mm {
\mathbb{M}U \ar@<1.5ex>[rr]^{\mathbb{M}_*i \phantom{aaa}} \ar@{}[rr]|-{\bot} && \mathbb{M}(U+V)\ar@<1.5ex>[ll]^{\mathbb{M}i \phantom{aaa}} \ar@<1.5ex>[rr]^{\phantom{aaa} \mathbb{M}j} \ar@{}[rr]|-{\top} && \mathbb{M}V \ar@<1.5ex>[ll]^{\phantom{aaa} \mathbb{M}_*j }
}
\end{eqnarray*}
exhibits a (bicategorical) direct sum in $\mathrm{Cat}_{+}$.
\end{corollary}

\begin{corollary}\label{aboutC_R}
For all $f : X'\to X$ in $\CC$, the functor $\mathbb{C}_{\CR}f : \mathbb{C}_{\CR}X \to  \mathbb{C}_{\CR}X'$ preserves finite coproducts and the pseudofunctor
\begin{eqnarray*}
 \CC^{\mathrm{op}}\times \CC\xra{(-)_*\times (-)_*}  (\mathrm{Spn}_{\CR})^{\mathrm{op}}\times \mathrm{Spn}_{\CR}\xra{\mathrm{Spn}_{\CR}(-,-)}\mathrm{Cat}
\end{eqnarray*}
factors through the inclusion $\mathrm{Cat}_{+}\hookrightarrow \mathrm{Cat}$ of the sub-2-category of categories with finite coproducts and functors which preserve them.
\end{corollary}

\begin{remark} Example~\ref{lextensiveE=Rcont} told of a hom enrichment of $\mathrm{Spn}_{\CR}\CE$ in $\mathrm{Cat}_{+}$ in the ordinary sense. For $\mathrm{Spn}_{\CR}\CC$,
we have an example of enrichment in a skew base which is an objective version of a base
observed in 2012 by 
Cockett and Lack in connection with left additive (= left linear over the ring of integers) categories; see page 1101 of \cite{GarLem}. 
What we notice is that the sub-2-category $\mathrm{Cat}_{+}$ of $\mathrm{Cat}$ has a skew closed structure \cite{116} with internal hom of $A,B\in \mathrm{Cat}_{+}$ taken as the ordinary functor category $[\CA,\CB]$ where the finite coproducts are pointwise;
write $\mathrm{Cat}_{+ \mathrm{sk}}$ for $\mathrm{Cat}_{+}$ equipped with this skew-monoidal structure. Then Proposition~\ref{fincoprodsglobloc} and Corollary~\ref{aboutC_R} tell us that $\mathrm{Spn}_{\CR}\CC$ is a skew
$\mathrm{Cat}_{+ \mathrm{sk}}$-enrichment of $\CC$ in the sense of Campbell \cite{Campb2018} where
the enriched hom is $\mathrm{Spn}_{\CR}\CC(X,Y)$. 

In particular, from Theorem~\ref{plycoprodthm}, we have that $\tensor[_{\CL}]{\mathrm{Ply}}{_\CR}\CE$ is a skew enrichment of $\tensor[_{\CL}]{\mathrm{Spn}}{}\CE$ in $\mathrm{Cat}_{+ \mathrm{sk}}$. 
\end{remark}

We end this section with a consequence of Proposition~\ref{Rcocomrestricts}.

\begin{proposition}\label{monadicityofobjMackey} The $\CR$-cocompletion pseudomonad $\overrightarrow{(-)}$ on $\mathrm{CAT}\CC$
restricts to a pseudomonad on the bicategory $\mathrm{fpCAT}\CC$ of
finite biproduct preserving pseudofunctors $\CC^{\mathrm{op}}\to \mathrm{Cat}$. The Eilenberg-Moore pseudo-algebras
are the objective $\CR$-Mackey functors.  
\end{proposition}

\section{A monoidal structure on $\mathrm{OMky}_{\CR}\CE$}\label{AmsoOMky}
The construction $[\mathbb{Y},\mathbb{Z}]_{\CR}$ defined and discussed in Proposition~\ref{inthom} and Example~\ref{inthomMky} provides a closed structure 
(in the sense of \cite{60}) on the bicategory $\mathrm{OMky}_{\CR}\CC$. 
The closed structure is symmetric monoidal: there is a symmetric tensor product
$\mathbb{M}\star\mathbb{N}$ defined by a pseudonatural equivalence
\begin{eqnarray*}
\mathrm{OMky}_{\CR}\CC (\mathbb{M}\star\mathbb{N}, \mathbb{L}) \simeq \mathrm{OMky}_{\CR}\CC (\mathbb{M}, [\mathbb{N},\mathbb{L}]_{\CR}) \ . 
\end{eqnarray*}
There is a bicategorically universal pseudonatural transformation 
$$\omega_{\mathbb{M},\mathbb{N}} : \mathbb{M}\times \mathbb{N} \to \mathbb{M}\star\mathbb{N}$$ 
which is $\CR$-cocontinuous in each variable separately.
The unit for this tensor product is the objective $\CR$-Burnside functor $\mathbb{C}_{\CR}$ (see Example~\eqref{oRBf}); for
\begin{eqnarray*}
[\mathbb{C}_{\CR},\mathbb{L}]_{\CR} & = & \mathrm{CAT}_{\CR}\CC(\overrightarrow{\mathbf{1}},\mathbb{L}(-\times W)) \\
& \simeq & \mathrm{CAT}\CC(\mathbf{1},\mathbb{L}(-\times W)) \\
& \simeq & \mathrm{lim} \ \mathbb{L}(-\times W) \\
& \simeq & \mathbb{L}(\mathbf{1}\times W) \\
& \simeq & \mathbb{L}W
\end{eqnarray*}
 where the (bicategorical) limit is over the category $\CC^{\mathrm{op}}$ which has the
(bicategorical) initial object $\mathbf{1}$. 
 
 The monoidales (pseudomonoids) for the convolution
monoidal structure $\mathbb{A}\star\mathbb{B}$ we call {\em objective $\CR$-Green functors on $\CC$} and write $\mathrm{OGrn}_{\CR}\CC$ 
for the bicategory they, with the strong monoidal morphisms, form.
An objective $\CR$-Green functor is {\em symmetric} when it is symmetric as a monoidale.

The tensor unit for a symmetric monoidal structure is always a symmetric monoidale. 
So the objective $\CR$-Burnside functor $\mathbb{C}_{\CR}$ is an example of an objective 
$\CR$-Green functor.    

Let us now describe the symmetric monoidal structure obtained by transporting the one on $\mathrm{OMky}_{\CR}\CC$ across the biequivalence \eqref{RspnclassCM}.  

By Lemma~\ref{finprodofmorphsinL}, we see that finite product in $\CC$ supplies the
bicategory $\mathrm{Spn}_{\CR}\CC$ with a symmetric monoidal structure. 
Using Proposition~\ref{Catsub+fact}, we see that the bicategory $\mathrm{fpCAT}\mathrm{Spn}_{\CR}\CC$
is the $\mathrm{Cat}_{+}$-enriched bicategory $\mathrm{Ps}((\mathrm{Spn}_{\CR}\CC)^{\mathrm{op}}, \mathrm{Cat}_{+})_{+}$ of such $\mathrm{Cat}_{+}$-enriched pseudofunctors. The desired transported symmetric closed monoidal structure
$\mathbb{A}\star\mathbb{B}$ on $\mathrm{fpCAT}\mathrm{Spn}_{\CR}\CC$ is the Day convolution structure with internal hom
$[\mathbb{B},\mathbb{D}]$ defined by
\begin{eqnarray*}
[\mathbb{B},\mathbb{D}]W
\ \simeq \ \mathrm{CAT}(\mathrm{Spn}_{\CR}\CC)(\mathbb{B},\mathbb{C}(-\times W))\ .
\end{eqnarray*}

\section{An objective Mackey functor of Mackey functors}\label{ToMfoMf}

Take a lextensive category $\CE$ equipped with a protocalibration $\CR$ closed under finite coproducts. We will consider the bicategory $\mathrm{Spn}_{\CR}{\CE}$ in the notation of Section~\ref{RSaRc} and will write $\mathrm{clSpn}_{\CR}{\CE}$ for the classifying category of $\mathrm{Spn}_{\CR}{\CE}$; 
the morphisms are isomorphism classes of spans (see \cite{Ben1967}). We write
$\mathrm{Spn}{\CE}$ and $\mathrm{clSpn}{\CE}$ for $\mathrm{clSpn}_{\CR}{\CE}$ and $\mathrm{clSpn}_{\CR}{\CE}$ in the case where $\CR = \CE$. 

As we have pointed out, the bicategory $\mathrm{Spn}_{\CR}{\CE}$ is enriched in the closed 
symmetric monoidal bicategory $\mathrm{Cat}_{+}$ of
categories with finite coproducts and finite-coproduct-preserving functors. The symmetric monoidal
structure on the bicategory $\mathrm{Spn}_{\CR}{\CE}$ is $\mathrm{Cat}_{+}$-enriched. The more abstract category $\mathrm{clSpn}_{\CR}{\CE}$, along with its symmetric 
monoidal structure, is enriched in the usual symmetric closed monoidal category $\mathrm{CMon}$ of commutative monoids.  

An {\em $\CR$-Mackey functor on $\CE$} is a $\mathrm{CMon}$-enriched (or equally, a finite-product-preserving) functor 
$$M : (\mathrm{clSpn}_{\CR}{\CE})^{\mathrm{op}}\to \mathrm{CMon} \ .$$
This is consistent with Lindner \cite{Lind} who looked at the case $\CR=\CE$.

\begin{example} \begin{itemize} 
\item[i.] Let $G$ be a finite group, let $\CE = [G,\mathrm{set}]$ be the category of finite $G$-sets, and let $R : G \to \mathrm{vect}$ be a finite-dimensional representation of $G$.
Define $M : \mathrm{clSpn}{\CE}\to \mathrm{CMon}$ by $MX = \mathrm{Set}^G(X,R)$ and, for $\tau : X\to R$ and $y\in Y$, 
$$M(X\xra{[u,S,v]}Y)(X\xra{\tau} R)(y) = \sum_{v(s) =y}{\tau(u(s))}\ .$$
\item[ii.] There is the {\em Burnside functor} $E_{\CR} : (\mathrm{clSpn}_{\CR}{\CE})^{\mathrm{op}}\to \mathrm{CMon}$ where 
$$E(X) = (\mathrm{clSpn}_{\CR}{\CE})(X,1) = \mathbb{E}_{\CR}(X)_{\cong}$$
is the set of isomorphism classes of the slice category $\CE_{/X}$ made into a commutative monoid with addition induced
by coproduct.  
\end{itemize}  
\end{example} 
Put
$$\mathrm{Mky}_{\CR}\CE = [(\mathrm{clSpn}_{\CR}{\CE})^{\mathrm{op}},\mathrm{CMon}]_+ \ ,$$
the $\mathrm{CMon}$-enriched functor category. 
It is a monoidal $\mathrm{CMon}$-enriched category by 
Day convolution \cite{DayConv}. The monoids for this convolution structure are {\em $\CR$-Green functors on $\CE$}; see \cite{93} for more details in the case $\CR = \CE$. We write $\mathrm{Grn}_{\CR}\CE$ for the $\mathrm{CMon}$-enriched category of Green
functors. The Burnside functor $E_{\CR}$ is an example. 

Let $\mathrm{RLxt}$ denote the 2-category whose objects are pairs $(\CA , \CO)$ consisting of
a lextensive category $\CA$ and a protocalibration $\CO$ closed under finite coproducts,
and, whose morphisms $F : (\CA , \CO)\to (\CA' , \CO')$ are functors $F : \CA\to \CA'$ preserving
pullback and finite coproducts, and taking $\CO$ to $\CO'$. Then we have a pseudofunctor $$\mathrm{Mky_{R}} : \mathrm{RLxt} \to \mathrm{Cat}_+$$ taking
$(\CA , \CO)$ to $\mathrm{Mky}_{\CO}\CA$ regarded as a category with finite coproducts. 

It is useful to notice that the objective Burnside functor $\mathbb{E}_{\CR} : \CE^{\mathrm{op}}\to \mathrm{Cat}_+$ lands in the category $\mathrm{RLxt}$ over $\mathrm{Cat}_{+}$ by equipping
each lextensive category $\mathbb{E}_{\CR}(X) (\subseteq \CE_{/X})$ with the protocalibration $\CR_{/X}$. 

\begin{proposition}
The composite pseudofunctor 
$$\mathrm{Mky}_{\mathrm{R}\CR} : = \CE^{\mathrm{op}}\xra{\mathbb{E}_{\CR}}\mathrm{RLxt}\xra{\mathrm{Mky_{R}}}\mathrm{Cat}_{+}$$
becomes a symmetric objective $\CR$-Green functor when equipped with the multiplication induced by the finite-coproduct-preserving-in-each-variable functors
$$\mathrm{Mky_R}\mathbb{E}_{\CR}X \times \mathrm{Mky_R}\mathbb{E}_{\CR}Y \to\mathrm{Mky_R}\mathbb{E}_{\CR}(X\times Y)$$
taking $(M,N)$ to $M\boxtimes N$ where $$(M\boxtimes N)(S\xra{(a,b)}X\times Y) = M(S\xra{a}X)\otimes N(S\xra{b}Y) \ .$$  
\end{proposition}

\begin{remark}\label{KvaluedobjMackey} In much of what we have just said the values of an objective Mackey functor could land
in any bicategory $\CK$ with homs enriched in $\mathrm{Cat}_{+}$ rather than in $\mathrm{Cat}_{+}$ itself.
That is, we can define a {\em $\CK$-valued objective Mackey functor on $\CE$} to be an $\CR$-cocomplete $\mathrm{Cat}_{+}$-enriched pseudofunctor  
$\mathbb{M} : \CE^{\mathrm{op}}\to\CK$. We write $\mathrm{OMky}_{\CR}(\CE,\CK)$ for the bicategory of these. 
For the convolution tensor product to exist on this bicategory,
we require $\CK$ to be appropriately cocomplete. 
We could do a similar thing in the abstract case by replacing $\mathrm{CMon}$
by a general $\mathrm{CMon}$-enriched category with finite coproducts.  
\end{remark}

\section{Objective Tambara functors}
Suppose $\CE$ is finitely complete and the pair $(\CL,\CR)$ satisfies the hypotheses of
Lemma~\ref{Rextensivity}.

\begin{definition}
An {\em objective $(\CL,\CR)$-Tambara functor} on $\CE$ is a 
finite biproduct preserving pseudofunctor $\mathbb{T} : \CE^{\mathrm{op}} \to \mathrm{Cat}$ which is in $\tensor[_{\CL}]{\mathrm{CAT}}{_\CR}\CE$.
We define the bicategory of objective $(\CL,\CR)$-Tambara functors by 
$$\tensor[_{\CL}]{\mathrm{Tmb}}{_\CR}\CE : = \mathrm{fp}\tensor[_{\CL}]{\mathrm{CAT}}{_\CR}\CE \ .$$ 
By Proposition~\ref{RspnclassC} and Theorem~\ref{PlyasSpnSpn}, 
\begin{eqnarray}\label{RspnclassCT}
\mathrm{fpCAT}\left( (\tensor[_{\CL}]{\mathrm{Ply}}{_\CR}\CE)^{\mathrm{op}}\right) \sim \tensor[_{\CL}]{\mathrm{Tmb}}{_\CR}\CE \ .
\end{eqnarray}
We write $\mathbb{T}_{\mathrm{pl}} : \mathrm{Ply}\CE \to \mathrm{Cat}$ 
for the inverse image of $\mathbb{T}$ under this biequivalence.  
\end{definition}
\begin{example} If $\mathbb{Z}$ is a $(\CL,\CR)$-Tambara functor then so is $[\mathbb{Y},\mathbb{Z}]_{\CR}$ for
all $\mathbb{Y}\in \mathrm{CAT}_{\CR}\CE$. 
\end{example}
\begin{example} Take $\CE$ to be a cartesian closed category $\mathrm{cat}$ of small categories and functors.
Let $\CX$ be a complete and cocomplete category. Then $\mathrm{Cat}(-,\CX) : \mathrm{cat}^{\mathrm{op}} \to \mathrm{Cat}$
is a $(\CP,\CE)$-Tambara functor.  
\end{example}
\begin{lemma} $\CE$ is $\CL$-extensive. 
\end{lemma}
\begin{proof} The inclusion $\mathbb{E}_{\CL}\to \mathbb{E}_{\CR}$ is componentwise fully
faithful and the $\CR$-extensivity equivalences restrict to $\mathbb{E}_{\CL}$.   
\end{proof}
Therefore Proposition~\ref{Rcocomrestricts} applies to both $\CR$ and $\CL$ yielding:

\begin{corollary} If $\mathbb{X}\in \mathrm{CAT}$ preserves biproducts then so do
$\overleftarrow{\mathbb{X}}$ and $\overrightarrow{\overleftarrow{\mathbb{X}}}$. 
\end{corollary}

\begin{proposition}\label{monadicityofobjTambara} The composite pseudomonad $\overrightarrow{\overleftarrow{(-)}}$ on $\mathrm{CAT}\CE$ (see Section~\ref{Adlfc})
restricts to a pseudomonad on the bicategory $\mathrm{fpCAT}\CE$ of
finite product preserving pseudofunctors $\CE^{\mathrm{op}}\to \mathrm{Cat}$. 
The pseudo-algebras are the objective $(\CL,\CR)$-Tambara functors. 
\end{proposition}

An {\em $(\CL,\CR)$-Tambara functor on $\CE$} is a finite-product-preserving functor 
$$T : \mathrm{cl}\tensor[_{\CL}]{\mathrm{Ply}}{_\CR}\CE\to \mathrm{CMon} \ .$$
This concept, for the case $\CL=\CR=\CE= G\text{-}\mathrm{set}$, arose in the paper \cite{Tamb1993} by Tambara; also see \cite{Maz, Hoyer, Ch} for a version of $\CL$ and $\CR$. 

\begin{center}
--------------------------------------------------------
\end{center}

\appendix

\end{document}